\documentclass[11pt]{article}
\usepackage{amsmath, upgreek}
\usepackage{amssymb}
\usepackage{amscd}
\usepackage{xspace}
\usepackage{verbatim}
\newcommand{\mysection}[1]{
\section{#1}\setcounter{equation}{0}}
\title{\bf Initial trace of positive solutions of a class of  degenerate heat equation with absorption}
\author{{\bf Tai Nguyen Phuoc}\qquad
 {\bf Laurent V\'eron}\\[2mm]
{\small Laboratoire de Math\'ematiques et Physique Th\'eorique, }\\
{\small  Universit\'e Fran\c{c}ois Rabelais,  Tours,  FRANCE}}

\date{}
\begin{document}
 \maketitle


\newcommand{\txt}[1]{\;\text{ #1 }\;}
\newcommand{\tbf}{\textbf}
\newcommand{\tit}{\textit}
\newcommand{\tsc}{\textsc}
\newcommand{\trm}{\textrm}
\newcommand{\mbf}{\mathbf}
\newcommand{\mrm}{\mathrm}
\newcommand{\bsym}{\boldsymbol}
\newcommand{\scs}{\scriptstyle}
\newcommand{\sss}{\scriptscriptstyle}
\newcommand{\txts}{\textstyle}
\newcommand{\dsps}{\displaystyle}
\newcommand{\fnz}{\footnotesize}
\newcommand{\scz}{\scriptsize}
\newcommand{\be}{\begin{equation}}
\newcommand{\bel}[1]{\begin{equation}\label{#1}}
\newcommand{\ee}{\end{equation}}
\newcommand{\eqnl}[2]{\begin{equation}\label{#1}{#2}\end{equation}}
\newcommand{\barr}{\begin{eqnarray}}
\newcommand{\earr}{\end{eqnarray}}
\newcommand{\bars}{\begin{eqnarray*}}
\newcommand{\ears}{\end{eqnarray*}}
\newcommand{\nnu}{\nonumber \\}
\newtheorem{subn}{\name}
\renewcommand{\thesubn}{}
\newcommand{\bsn}[1]{\def\name{#1}\begin{subn}}
\newcommand{\esn}{\end{subn}}
\newtheorem{sub}{\name}[section]
\newcommand{\dn}[1]{\def\name{#1}}   
\newcommand{\bs}{\begin{sub}}
\newcommand{\es}{\end{sub}}
\newcommand{\bsl}[1]{\begin{sub}\label{#1}}
\newcommand{\bth}[1]{\def\name{Theorem}
\begin{sub}\label{t:#1}}
\newcommand{\blemma}[1]{\def\name{Lemma}
\begin{sub}\label{l:#1}}
\newcommand{\bcor}[1]{\def\name{Corollary}
\begin{sub}\label{c:#1}}
\newcommand{\bdef}[1]{\def\name{Definition}
\begin{sub}\label{d:#1}}
\newcommand{\bprop}[1]{\def\name{Proposition}
\begin{sub}\label{p:#1}}
\newcommand{\R}{\eqref}
\newcommand{\rth}[1]{Theorem~\ref{t:#1}}
\newcommand{\rlemma}[1]{Lemma~\ref{l:#1}}
\newcommand{\rcor}[1]{Corollary~\ref{c:#1}}
\newcommand{\rdef}[1]{Definition~\ref{d:#1}}
\newcommand{\rprop}[1]{Proposition~\ref{p:#1}}
\newcommand{\BA}{\begin{array}}
\newcommand{\EA}{\end{array}}
\newcommand{\BAN}{\renewcommand{\arraystretch}{1.2}
\setlength{\arraycolsep}{2pt}\begin{array}}
\newcommand{\BAV}[2]{\renewcommand{\arraystretch}{#1}
\setlength{\arraycolsep}{#2}\begin{array}}
\newcommand{\BSA}{\begin{subarray}}
\newcommand{\ESA}{\end{subarray}}
\newcommand{\BAL}{\begin{aligned}}
\newcommand{\EAL}{\end{aligned}}
\newcommand{\BALG}{\begin{alignat}}
\newcommand{\EALG}{\end{alignat}}
\newcommand{\BALGN}{\begin{alignat*}}
\newcommand{\EALGN}{\end{alignat*}}
\newcommand{\note}[1]{\textit{#1.}\hspace{2mm}}
\newcommand{\Proof}{\note{Proof}}
\newcommand{\qeda}{\hspace{10mm}\hfill $\square$}
\newcommand{\qed}{\\
${}$ \hfill $\square$}
\newcommand{\Remark}{\note{Remark}}
\newcommand{\modin}{$\,$\\[-4mm] \indent}
\newcommand{\forevery}{\quad \forall}
\newcommand{\set}[1]{\{#1\}}
\newcommand{\setdef}[2]{\{\,#1:\,#2\,\}}
\newcommand{\setm}[2]{\{\,#1\mid #2\,\}}
\newcommand{\mt}{\mapsto}
\newcommand{\lra}{\longrightarrow}
\newcommand{\lla}{\longleftarrow}
\newcommand{\llra}{\longleftrightarrow}
\newcommand{\Lra}{\Longrightarrow}
\newcommand{\Lla}{\Longleftarrow}
\newcommand{\Llra}{\Longleftrightarrow}
\newcommand{\warrow}{\rightharpoonup}
\newcommand{
\paran}[1]{\left (#1 \right )}
\newcommand{\sqbr}[1]{\left [#1 \right ]}
\newcommand{\curlybr}[1]{\left \{#1 \right \}}
\newcommand{\abs}[1]{\left |#1\right |}
\newcommand{\norm}[1]{\left \|#1\right \|}
\newcommand{
\paranb}[1]{\big (#1 \big )}
\newcommand{\lsqbrb}[1]{\big [#1 \big ]}
\newcommand{\lcurlybrb}[1]{\big \{#1 \big \}}
\newcommand{\absb}[1]{\big |#1\big |}
\newcommand{\normb}[1]{\big \|#1\big \|}
\newcommand{
\paranB}[1]{\Big (#1 \Big )}
\newcommand{\absB}[1]{\Big |#1\Big |}
\newcommand{\normB}[1]{\Big \|#1\Big \|}
\newcommand{\produal}[1]{\langle #1 \rangle}

\newcommand{\thkl}{\rule[-.5mm]{.3mm}{3mm}}
\newcommand{\thknorm}[1]{\thkl #1 \thkl\,}
\newcommand{\trinorm}[1]{|\!|\!| #1 |\!|\!|\,}
\newcommand{\bang}[1]{\langle #1 \rangle}
\def\angb<#1>{\langle #1 \rangle}
\newcommand{\vstrut}[1]{\rule{0mm}{#1}}
\newcommand{\rec}[1]{\frac{1}{#1}}
\newcommand{\opname}[1]{\mbox{\rm #1}\,}
\newcommand{\supp}{\opname{supp}}
\newcommand{\dist}{\opname{dist}}
\newcommand{\myfrac}[2]{{\displaystyle \frac{#1}{#2} }}
\newcommand{\myint}[2]{{\displaystyle \int_{#1}^{#2}}}
\newcommand{\Lim}[1]{{\displaystyle \lim_{#1}}}
\newcommand{\mysum}[2]{{\displaystyle \sum_{#1}^{#2}}}
\newcommand {\dint}{{\displaystyle \myint\!\!\myint}}
\newcommand{\q}{\quad}
\newcommand{\qq}{\qquad}
\newcommand{\hsp}[1]{\hspace{#1mm}}
\newcommand{\vsp}[1]{\vspace{#1mm}}
\newcommand{\ity}{\infty}
\newcommand{\prt}{\partial}
\newcommand{\sms}{\setminus}
\newcommand{\ems}{\emptyset}
\newcommand{\ti}{\times}
\newcommand{\pr}{^\prime}
\newcommand{\ppr}{^{\prime\prime}}
\newcommand{\tl}{\tilde}
\newcommand{\sbs}{\subset}
\newcommand{\sbeq}{\subseteq}
\newcommand{\nind}{\noindent}
\newcommand{\ind}{\indent}
\newcommand{\ovl}{\overline}
\newcommand{\unl}{\underline}
\newcommand{\nin}{\not\in}
\newcommand{\pfrac}[2]{\genfrac{(}{)}{}{}{#1}{#2}}

\def\ga{\alpha}     \def\gb{\beta}       \def\gg{\gamma}
\def\gc{\chi}       \def\gd{\delta}      \def\ge{\epsilon}
\def\gth{\theta}                         \def\vge{\varepsilon}
\def\gf{\phi}       \def\vgf{\varphi}    \def\gh{\eta}
\def\gi{\iota}      \def\gk{\kappa}      \def\gl{\lambda}
\def\gm{\mu}        \def\gn{\nu}         \def\gp{\pi}
\def\vgp{\varpi}    \def\gr{\rho}        \def\vgr{\varrho}
\def\gs{\sigma}     \def\vgs{\varsigma}  \def\gt{\tau}
\def\gu{\upsilon}   \def\gv{\vartheta}   \def\gw{\omega}
\def\gx{\xi}        \def\gy{\psi}        \def\gz{\zeta}
\def\Gg{\Gamma}     \def\Gd{\Delta}      \def\Gf{\Phi}
\def\Gth{\Theta}
\def\Gl{\Lambda}    \def\Gs{\Sigma}      \def\Gp{\Pi}
\def\Gw{\Omega}     \def\Gx{\Xi}         \def\Gy{\Psi}

\def\CS{{\mathcal S}}   \def\CM{{\mathcal M}}   \def\CN{{\mathcal N}}
\def\CR{{\mathcal R}}   \def\CO{{\mathcal O}}   \def\CP{{\mathcal P}}
\def\CA{{\mathcal A}}   \def\CB{{\mathcal B}}   \def\CC{{\mathcal C}}
\def\CD{{\mathcal D}}   \def\CE{{\mathcal E}}   \def\CF{{\mathcal F}}
\def\CG{{\mathcal G}}   \def\CH{{\mathcal H}}   \def\CI{{\mathcal I}}
\def\CJ{{\mathcal J}}   \def\CK{{\mathcal K}}   \def\CL{{\mathcal L}}
\def\CT{{\mathcal T}}   \def\CU{{\mathcal U}}   \def\CV{{\mathcal V}}
\def\CZ{{\mathcal Z}}   \def\CX{{\mathcal X}}   \def\CY{{\mathcal Y}}
\def\CW{{\mathcal W}} \def\CQ{{\mathcal Q}}
\def\BBA {\mathbb A}   \def\BBb {\mathbb B}    \def\BBC {\mathbb C}
\def\BBD {\mathbb D}   \def\BBE {\mathbb E}    \def\BBF {\mathbb F}
\def\BBG {\mathbb G}   \def\BBH {\mathbb H}    \def\BBI {\mathbb I}
\def\BBJ {\mathbb J}   \def\BBK {\mathbb K}    \def\BBL {\mathbb L}
\def\BBM {\mathbb M}   \def\BBN {\mathbb N}    \def\BBO {\mathbb O}
\def\BBP {\mathbb P}   \def\BBR {\mathbb R}    \def\BBS {\mathbb S}
\def\BBT {\mathbb T}   \def\BBU {\mathbb U}    \def\BBV {\mathbb V}
\def\BBW {\mathbb W}   \def\BBX {\mathbb X}    \def\BBY {\mathbb Y}
\def\BBZ {\mathbb Z}

\def\GTA {\mathfrak A}   \def\GTB {\mathfrak B}    \def\GTC {\mathfrak C}
\def\GTD {\mathfrak D}   \def\GTE {\mathfrak E}    \def\GTF {\mathfrak F}
\def\GTG {\mathfrak G}   \def\GTH {\mathfrak H}    \def\GTI {\mathfrak I}
\def\GTJ {\mathfrak J}   \def\GTK {\mathfrak K}    \def\GTL {\mathfrak L}
\def\GTM {\mathfrak M}   \def\GTN {\mathfrak N}    \def\GTO {\mathfrak O}
\def\GTP {\mathfrak P}   \def\GTR {\mathfrak R}    \def\GTS {\mathfrak S}
\def\GTT {\mathfrak T}   \def\GTU {\mathfrak U}    \def\GTV {\mathfrak V}
\def\GTW {\mathfrak W}   \def\GTX {\mathfrak X}    \def\GTY {\mathfrak Y}
\def\GTZ {\mathfrak Z}   \def\GTQ {\mathfrak Q}

\font\Sym= msam10 
\def\SYM#1{\hbox{\Sym #1}}
\newcommand{\bdw}{\prt\Gw\xspace}
\tableofcontents
\medskip
\begin{abstract}
We study the  initial value problem with unbounded nonnegative functions or measures  for the equation $ \prt_tu-\Gd_p u+f(u)=0$ in $\BBR^N\ti(0,\infty)$ where $p>1$, $\Gd_p u = \text{div}(\abs {\nabla u}^{p-2} \nabla u )$ and $f$ is a continuous, nondecreasing nonnegative function such that $f(0)=0$. In the case $p>\frac{2N}{N+1}$, we provide a sufficient condition on $f$ for existence and uniqueness of the solutions satisfying the initial data $k\gd_0$ and we study their limit when $k\to\infty$, according $f^{-1}$ and $F^{-1/p}$ are integrable or not at infinity, where $F(s)=\int_0^s f(\gs)d\gs$. We also give new results dealing with non uniqueness for the initial value problem with unbounded initial data. If $p>2$, we prove that, for a large class of nonlinearities $f$, any positive solution admits an initial trace in the class of positive Borel measures. As a model case we consider the case $f(u)=u^\ga \ln^\gb(u+1)$, where $\ga>0$ and $\gb\geq 0$.
\end{abstract}
\mysection{Introduction}
The aim of this article is to study some qualitative properties of the positive solutions of 
\bel{A1} \prt_t u - \Gd_p u +f(u) = 0 \ee
in $Q_\ity:={\BBR}^N \ti (0,\ity)$ where $p>1$, $\Gd_p u = \text{div}(\abs {\nabla u}^{p-2} \nabla u )$ and $f$ is a continuous, nondecreasing function such that $f(0)=0=f^{-1}(0)$. The properties we are interested in are mainly: (a) the existence of fundamental solutions i.e. solutions with $k\gd_0$ as initial data and the behaviour of these solutions when $k\to\infty$; (b) the existence of an initial trace and its properties; (c) uniqueness and non-uniqueness results for the Cauchy problem. This type of questions have been considered in a previous paper of the authors \cite {NgV} in the semilinear case $p=2$. The breadcrumbs of this study lies in the existence of two types of specific solutions of $(\ref{A1})$. The first ones are the solutions $\gf:=\gf_a$ of the ODE
\bel{I-1} \phi' +f(\phi) = 0 \ee
defined on $[0,\infty)$ and subject to $\phi(0)=a\geq 0$; it is given by
\bel{I-2}
\myint{\phi(t)}{a}\myfrac{ds}{f(s)}.
\ee
The second ones are the solutions of the elliptic equation
\bel{I-3} -\Gd_p w +f(w) = 0, \ee
defined in $\BBR^N$ or in $\BBR^N\setminus\{0\}$.
It is well-known that the structure of the set of solutions of $(\ref{I-1})$ depends whether the following quantity
\bel{I-4}
J:=\myint{1}{\infty}\myfrac{ds}{f(s)}
\ee 
is finite or infinite. If $J<\infty$ there exists a maximal solution $\phi_\infty$ to $(\ref{I-1})$ defined on $(0,\infty)$ while no such solution exists if $J=\infty$ since $\lim_{a\to\infty}\gf_a(t)=\infty$. This maximal solution plays an important role since, by the maximum principle, it dominates any solution $u$ of $(\ref{A1})$ which satisfies
\bel{I-C}
\lim_{|x|\to\infty}u(x,t)=0
\ee
for all $t>0$, locally uniformly on $(0,\infty)$. Concerning $(\ref{I-3})$ we associate the quantity
\bel{I-5}
K:=\myint{1}{\infty}\myfrac{ds}{{F(s)}^{1/p}}.
\ee 
It is a consequence of the V\'azquez's extension of the Keller-Osserman condition (see \cite{Va}, \cite{Ke}) that if $K<\infty$, equation $(\ref{I-3})$ admits a maximal solution $W_{\BBR^N_*}$ in $\BBR^N\setminus\{0\}$. This solution is constructed as the limit, when $R\to\infty$ and $\ge\to 0$ of the solution $W:=W_{\ge,R}$ of $(\ref{I-3})$ in $\Gg_{\ge,R}:=B_R\setminus \overline B_\ge$, subject to the conditions 
$\lim_{|x|\downarrow\ge}W_{\ge,R}(x)=\infty$ and $\lim_{|x|\uparrow R}W_{\ge,R}(x)=\infty.$ 
On the contrary, if $K=\infty$, such functions $W_{\ge,R}$ and $W_{\BBR^N_*}$ do not exist, a situation which will be exploited in Section 3 for proving existence of global solutions of $(\ref{I-3})$ in $\BBR^N$. An additional natural growth assumption of $f$ that will be often made is the {\it super-additivity}
\bel{addit} f(s+s')\geq f(s)+f(s')\qquad\forall s,s'\geq 0, \ee
which, combined with the monotonicity of $f$, implies a minimal linear growth at infinity
\bel{superlin}\liminf_{s\to\infty}\myfrac{f(s)}{s}>0.\ee
If $p\geq 2$, $K<\infty$ jointly with $(\ref{addit})$ implies $J<\infty$, but this does not hold when $1<p<2$. When $p>2$ and $f$ satisfies $J<\ity$ and $K<\ity$, Kamin and V\'azquez proved universal estimates for solutions which vanish on $\BBR^N\ti\{0\}\setminus \{(0,0)\}$ (see \cite{KV2}). By a slight modification of the proof in \cite[Proposition 2.3 and Proposition 2.6]{NgV}, it is possible to extend their result to the case $p>1$. \qeda \medskip




\noindent{\bf Proposition} (Universal estimates) {\it Assume $p>1$ and $f$ satisfies $K<\infty$. Let $u\in C(\overline {Q_\infty}\setminus \{(0,0)\})$ be a solution of $(\ref{A1})$ in $Q_\infty$, which vanishes on 
$\BBR^N\ti\{0\}\setminus \{(0,0)\}$.Then
\bel{I-5-1}
u(x,t)\leq W_{\BBR^N_*}(x)\qquad\forall (x,t)\in Q_\infty.
\ee
 If we suppose moreover $J<\infty$ and that $(\ref{addit})$ holds, then}  
\bel{I-5-2}
u(x,t)\leq \min\left\{\gf_\infty(t),W_{\BBR^N_*}(x)\right\}\qquad\forall (x,t)\in Q_\infty.
\ee

When $K=\infty$, no such estimate exists since the function $w_a$ solution of $(\ref{I-7})$ is a stationnary solution of $(\ref{A1})$ with unbounded initial data. \medskip

In Section 2 we study the existence of the fundamental solutions $u_k$ and their behaviour when $k\to\infty$. Kamin and V\'azquez proved in \cite[Lemma 2.3 and Lemma 2.4]{KV2}, that if $p>2$ and 
\bel{CFS} \myint{1}{\ity}s^{-p-\frac{p}{N}}f(s)ds < \ity, \ee
then for any $k>0$, there exists a unique positive solution $u:=u_k$ to problem
	\bel{I-9} \left\{ \BA{lll} \prt_t u - \Delta_p u + f(u) &= 0 \qq &\text{in } Q_\ity \\ [2mm] \phantom{-----,}
	u(.,0) &= k\delta_0 & \text{in } \BBR^N. \EA \right. 
\ee
Furthermore the mapping $k\mapsto u_k$ is increasing.  Their existence proof heavily relies on the fact that, if we denote by $v:=v_k$ the fundamental (or Barenblattt-Prattle) solution of 
\bel{I-10}\left\{ \BA{lll} \prt_t v - \Delta_p v &= 0 \qq &\text{in } Q_\ity \\ [2mm] \phantom{--,}
	v(.,0) &= k\delta_0 & \text{in } \BBR^N, \EA \right.
\ee
then $v_k(.,t)$ is compactly supported in some ball $B_{\gd_k(t)}$, where $\gd_k(t)$ is explicit. Since $v_k$ is a natural supersolution for $(\ref{I-9})$, condition $(\ref{CFS})$ states that $f(v_k)\in L_{loc}^1(\overline{Q_\infty})$. When $2N/(N+1)<p\leq 2$, $v_k(x,t)>0$ for all $(x,t)\in Q_\infty$. It is already proved in \cite {MV} that, when $p=2$, condition $(\ref{CFS})$ yields to $f(v_k)\in L^1(Q_T)$. We prove here that this result also holds when $2N/(N+1)<p\leq 2$ and more precisely,
\bth {Fund}Assume $p>\frac{2N}{N+1}$ and $f$ satisfies $(\ref{CFS})$. Then there exists a unique positive solution $u:=u_k$ to problem $(\ref{I-9})$.
\es

In view of this result and the {\it a priori} estimates $(\ref{I-5-1})$ and $(\ref{I-5-2})$, it is natural to study the limit of $u_k$ when $k\to\infty$. We denote by $\CU_0$ the set of positive  $u\in C(\overline{Q_\infty}\setminus\{(0,0)\})$ which are solutions of $(\ref{A1})$ in $Q_\infty$, vanishes on the set $\{(x,0): x\neq 0\}$ and satisfies
	$$ \lim_{t \to 0}\myint{B_\ge}{}u(x,t)dx = \ity \forevery \ge>0. $$

\bth {Strong} Assume $p>2N/(N+1)$, $J<\infty$, $K<\infty$ and $(\ref{CFS})$ holds. Then ${\underline U=\dsps \lim_{k\to\infty}u_k}$ exists and it is the smallest element of $\CU_0$.
\es

When one, at least, of the above properties on $J$ and $K$ fails, the situation is much more complicated and fairly well understood only in the case where $f$ has a power-like or a logarithmic-power-like growth. We first note that \smallskip

\noindent (A) If $f(s)\sim s^\ga$ ($\ga>0$), then $J<\infty$ if and only if $\ga>1$, while $K<\infty$ if and only if $\ga>p-1$. Moreover $(\ref{CFS})$ holds if and only if $\ga<p(1+\frac{1}{N})-1$. \smallskip

\noindent (B) If $f(s)\sim s^\ga\ln^\gb (s+1)$ ($\ga,\gb>0$), then $J<\infty$ if and only if $\ga>1$ and $\gb>0$, or $\ga=1$ and $\gb>1$ while $K<\infty$ if and only if $\ga>p-1$ and $\gb>0$, or $\ga=p-1$ and $\gb>p$. Moreover $(\ref{CFS})$ holds if and only if $\ga<p(1+\frac{1}{N})-1$ and $\gb>0$.\medskip

\bth{IS1} Assume $p>2$ and $f(s)=s^\ga\ln^\gb(s+1)$ where $\ga \in (1,p-1)$ and $\gb>0$. Let $u_k$ be the solution of $(\ref{I-9})$. Then ${\dsps \lim_{k \to \ity}{u_k(x,t)}}=\gf_\ity(t)$ for every $(x,t) \in Q_\ity$.
\es

When $\ga=1$ the following phenomenon occurs.

\bth{IS2} Assume $p>2$ and $f(s)=s\ln^\gb(s+1)$ with $\gb>0$. Let $u_k$ be the solution of $(\ref{I-9})$. Then \smallskip

\noindent (i) If $\gb>1$ then ${\dsps \lim_{k \to \ity}{u_k(x,t)}}=\gf_\ity(t)$ for every $(x,t) \in Q_\ity$, \smallskip

\noindent (ii) If $0<\gb\leq 1$ then ${\dsps \lim_{k \to \ity}{u_k(x,t)}}=\ity$ for every $(x,t) \in Q_\ity$.
\es

 
Section 3 is devoted to study non-uniqueness of solutions of $(\ref{A1})$ with unbounded initial data. The starting observation is the following global existence result for solutions of $(\ref{I-3})$: 
\bth{glob} Assume $p>1$, $f$ is locally Lipschitz continuous and $K=\infty$. Then for any $a>0$, there exists a unique solution $w:=w_a$ to the problem
\bel{I-6}
-(r^{N-1}|w_r|^{p-2}w_r)_r+r^{N-1}f(w)=0\\
\ee 
defined on $[0,\infty)$ and satisfying $w(0)=a$, $w_r(0)=0$. It is given by
\bel{I-7}
w_a(r)=a+\myint{0}{r}H_p\left(s^{1-N}\myint{0}{s}\gt^{N-1}f(w_a(\gt))d\gt\right)ds
\ee 
where $H_p$ is the inverse function of $t\mapsto |t|^{p-2}t$.
\es
  This result extends to the general case $p>1$ a previous theorem of V\'azquez and V\'eron \cite{VaVe1} obtained  in the case $p=2$. 
The next theorem extends to the case $p\neq 2$ a previous result of the authors in the case $p=2$.
\bth{non-unique} Assume $p>2N/(N+1)$, $f$ is locally Lipschitz continuous, $J<\infty$ and $K=\infty$. For any function $u_0 \in C(Q_\ity)$ which satisfies
	\be w_a(\abs x)\leq u_0(x) \leq w_b(\abs x) \q \forall x \in {\BBR}^N	\ee
for some $0<a<b$, there exist at least two solutions $\unl u, \ovl u \in C(\ovl {Q_\ity})$ of $(\ref{A1})$ with initial value $u_0$. They satisfy respectively
	$$ 0\leq \unl u(x,t) \leq \min\{w_b(\abs x),\gf_\ity(t)\} \forevery (x,t) \in Q_\ity, $$
thus ${\dsps \lim_{t \to \ity}\unl u(x,t)=0}$, uniformly with respect to $x \in \BBR^N$, and 
	$$ w_a(\abs x) \leq \ovl u(x,t) \leq w_b(\abs x) \forevery (x,t) \in Q_\ity $$
thus ${\dsps \lim_{\abs x \to \ity} \ovl u(x,t)=\ity}$, uniformly with respect to $t \geq 0$.
\es 

In section 4 we prove an existence and stability result for the initial value problem
	\bel{I-4-1} \left\{ \BA{lll} \prt_t u - \Delta_p u + f(u) &= 0 \qq &\text{in } Q_\ity \\ [2mm] \phantom{-----,}
	u(.,0) &= \gm & \text{in } \BBR^N \EA \right. 
\ee
where $\gm \in \GTM_+^b(\BBR^N)$, the set of positive and bounded Radon measures in $\BBR^N$.

\bth{stab} Assume $p>\frac{2N}{N+1}$ and $f$ satisfies $(\ref{CFS})$. Then for any $\gm \in \GTM_+^b(\BBR^N)$ the problem 
$(\ref{I-4-1})$ admits a weak solution $u_{\gm}$. Moreover, if $\{\gm_n\}$ is a sequence of functions in $L_+^1(\BBR^N)$ with compact support, which converges to $\gm \in \GTM_+^b(\BBR^N)$ in the weak sense of measures, then the corresponding solutions $\{u_{\gm_n}\}$ of $(\ref{I-4-1})$ with initial data $\gm_n$ converge to some solution $u_\gm$ of $(\ref{I-4-1})$, strongly in $L_{loc}^1(\overline{Q_T})$ and locally uniformly in $Q_T:=\BBR^N \ti (0,T)$. Furthermore $\{f(u_{\gm_n})\}$ converges strongly to $f(u_\gm)$ in $L_{loc}^1(\overline{Q_T})$.
\es  

In Section 5, we discuss the {\it initial trace} of positive {\it weak solution} of $(\ref{A1})$. The power case $f(u)=u^q$ with $q>0$ was investigated by Bidaut-V\'eron, Chasseigne and V\'eron in \cite{BvCV}. They proved the existence of an initial trace in the class of positive Borel measures according to the different values of $p-1$ and $q$. Accordingly they studied the corresponding Cauchy problem with a given Borel measure as initial data. However their method was strongly based upon the fact that the nonlinearity was a power, which enabled to use H\"older inequality in order to show the domination of the absorption term over the other terms. In the present paper, we combine the ideas in \cite{BvCV} and \cite{NgV} with a  stability result for the Cauchy problem  and Harnack's inequality in the form of \cite{Di} to establish the following dichotomy result which is new even in the case $p=2$.

\bth{dichotomy} Assume $p\geq 2$ and $(\ref{CFS})$ holds. Let $u \in C(Q_T)$ be a positive weak solution of $(\ref{A1})$ in $Q_T$. Then for any $y \in \BBR^N$ the following alternative holds \medskip

\noindent{(i) } either 
\bel{di1}u(x,t)\geq \lim_{k\to\infty}u_k(x-y,t)\qquad\forall (x,t)\in Q_T,
\ee
\noindent{(ii) } or there exist an open neighborhood $U$ of $y$ and a Radon measure $\gm_U \in \GTM_+(U)$ such that
	\bel{di2} \lim_{t \to 0}\myint{U}{}u(x,t)\zeta(x)dx=\myint{U}{}\zeta d\gm_U \forevery \zeta \in C_c(U). \ee

\es 

Actually, since $(\ref{CFS})$ is verified, $(\ref{di1})$ is equivalent to the fact that, for any open neighborhood $U$ of $y$, there holds
\bel{di3}\limsup_{t\to 0}\myint{U}{}u(x,t)dx=\infty.
\ee
However, if $(\ref{CFS})$ is not verified, there only holds $(\ref{di1})\Longrightarrow (\ref{di3})$. 

The set of points $y$ such that $(\ref{di2})$ (resp. $(\ref{di3})$) holds is clearly open (resp. closed) and denoted by $\CR(u)$ (resp ($\CS(u)$). Using a partition of unity, there exists a unique Radon measure $\gm\in \GTM_+(\CR(u))$ such that 
\bel{di4} 
\lim_{t \to 0}\myint{\CR(u)}{}u(x,t)\zeta(x)dx=\myint{\CR(u)}{}\zeta d\gm \forevery \zeta \in C_c(\CR(u)). 
\ee
Owing to the above result we  define the {\it initial trace} of a positive solution $u$ $(\ref{A1})$ in $Q_T$ as the couple $(\CS(u),\gm)$ for which $(\ref{di2})$ and $(\ref{di3})$ holds and we denote it by  $tr_{_{\BBR^N}}(u)$. The set $\CS(u)$ is the {\it set of singular points} of $tr_{_{\BBR^N}}(u)$, while $\gm$ is the {\it regular part} of $tr_{_{\BBR^N}}(u)$. It is classical that any $\gn\in\mathfrak B^{reg}(\BBR^N)$, the set of positive outer regular Borel measures in $\BBR^N$, can be represented by a couple $(\CS,\gm)$ where $\CS$ is a closed subset of $\BBR^N$ and $\gm\in \GTM_+(\CR)$, where $\CR=\BBR^N\setminus\CS$, in the following way 
	$$ \nu(A)=\left\{ \BA{ll} \ity \qq & \text{if } A \cap \CS \ne \emptyset, \\
	\gm(A) & \text{if } A \sbs \CR, \EA \right. \forevery A \text{ Borel}. $$
Therefore \rth{dichotomy} means that $tr_{_{\BBR^N}}(u)\in \mathfrak B^{reg}(\BBR^N)$.

The initial trace can be made more precise when the Keller-Osserman-V\'azquez condition does not hold, and if we know  whether ${\dsps \lim_{k \to \ity}u_k}$ is equal to $\gf_\ity$ or is infinite.

\bth{tr+J-fin} Assume $p>2$ and $(\ref{CFS})$ holds and $u$ is a positive solution of $(\ref{A1})$ in $Q_{\infty}$. \smallskip

 \noindent I- If $J<\infty$ and $K=\infty$ are verified and ${\dsps \lim_{k \to \ity}u_k}=\gf_\ity$. Then either $tr_{_{\BBR^N}}(u)$ is  the Borel measure infinity $\gn_\infty$ which satisfies $\gn_\infty(\CO)=\infty$ for any non-empty open subset $\CO\subset \BBR^N$, or is a positive Radon measure $\gm$ on $\BBR^N$. \smallskip

 \noindent II- If $J=\infty$ and $K=\infty$ are verified and ${\dsps \lim_{k \to \ity}u_k}=\ity$. Then  $tr_{_{\BBR^N}}(u)$ is   a positive Radon measure $\gm$ on $\BBR^N$
\es

As a consequence of $I$, there exist infinitely many positive solutions $u$ of $(\ref{A1})$ in $Q_{\infty}$ such that $tr_{_{\BBR^N}}(u)=\gn_\infty$. By \rth{IS1}, \rth{IS2}, the previous results apply in particular if $f(s)=s^\ga\ln^\gb(s+1)$. 
\mysection{Isolated singularities}
Throughout the article $c_i$ denote positive constants depending on $N$, $p$, $f$ and sometimes other quantities such as test functions or particular exponents, the value of which may change from one occurrence to another.
\subsection{The semigroup approach}
We refer to \cite[p 117]{HeVa} for the detail of the Banach space framework for the construction of solutions of $(\ref{A1})$ in $Q_\infty$ with initial data in $L^1(\BBR^N)\cap L^\infty(\BBR^N)$. We set 
\bel{X1}
J(u)=\myint{\BBR^N}{}\left(\frac{1}{p}\abs{\nabla u}^p+F(u)\right)dx
\ee
when $u$ belongs to the domain $D(J)$ of $J$ which is the set of $u\in L^2(\BBR^N)$ such that $\nabla u\in L^p(\BBR^N)$ and $F(u)\in L^1(\BBR^N)$, and $J(u)=\infty$ if $u\notin D(J)$. Then $J$ is a proper convex lower semicontinuous function in $L^2(\BBR^N)$. Its sub-differential $A$ is defined by its domain $D(A)$ which is the set 
of $u\in L^2(\BBR^N)$ such that $\nabla u\in L^p(\BBR^N)$ and $F(u)\in L^1(\BBR^N)$ with the property that $-\Gd_pu+f(u)\in L^2(\BBR^N)$ and 
\begin{equation}\label{X2}
-\myint{\BBR^N}{}v\Gd_pu dx=\myint{\BBR^N}{}\abs{\nabla u}^{p-2}\nabla u.\nabla vdx\qquad\forall v\in D(J),
\end{equation}
and by its expression
\bel{X3}
Au=-\Gd_pu+f(u)\qquad\forall u\in D(A).
\ee
Notice that $(\ref{X2})$ implies that $vf(u)\in L^1(\BBR^N)$ for all $v\in D(J)$. The restriction of the operator $A$ is accretive in $L^1(\BBR^N)$ and in $L^\infty(\BBR^N)$, hence in every $L^q(\BBR^N)$. The operator $A_q$ defined in $L^q(\BBR^N)$  is the closure in $L^q(\BBR^N)$ of the restriction of $A$ to $L^q(\BBR^N)$. It is a m-accretive operator, with domain $D(A_q)$. Since $C^\infty_0(\BBR^N)\subset D(A_q)$, $D(A_q)$ is dense in $L^q(\BBR^N)$.  If $u_0\in L^q$ the generalized solution $u$ to
\bel{X4}\left\{\BA {l}
\myfrac{du}{dt}+A_qu=0\qquad\text{in }(0,\infty)\\
\phantom{\frac{du}{dt}-.}
u(0)=u_0
\EA
\right.\ee
is obtained by the Crandall-Liggett scheme
\bel{X4*}\BA {l}
\myfrac{u_i-u_{i-1}}{h}+A_qu_i=0\qquad\text{in }i=0,1,...
\EA\ee
when we let $h\to 0$, in the sense that the continuous piecewise linear function $U_h$ defined by $U_h(ih)=u_i$ converges to $u$ in the $C([0,T],L^q(\BBR^N))$-topology, for every $T>0$. Furthermore, if $q=2$ and $u_0\in D(A_2)$ (resp.  $u_0\in L^2(\BBR^N)$), then $\frac{dU_h}{dt}$ converges to $\frac{du}{dt}$ in $L^2([0,T],L^2(\BBR^N))$ (resp. $L^2([0,T],L^2(\BBR^N);tdt)$), see \cite{Ve1}. We shall denote by $\{S^{A_q}(t)\}_{t>0}$ the semigroup of contractions of $L^q(\BBR^N)$ generated by $-A_q$ thru the Crandall-Liggett Theorem \cite{CL}.\smallskip

 An important property \cite [Lemma 2]{HeVa} is that if $w\in L^1(\BBR^N)$  satisfies
\bel{X5}\BA {l}
A_1w+\gs w=h
\EA\ee
where $\gs>0$ and $h\in L^1(\BBR^N)$, then 
\bel{X6}\BA {l}
\myint{\BBR^N}{}A_1 wdx=0.
\EA\ee

\bdef{semigroup} (i) A function $u\in C([\gd,\infty);L^1(\BBR^N))$  where $\gd\geq 0$ is a semigroup solution $(\ref{A1})$ on $(\gd,\infty)$ if for any $t\geq \gd$ there holds $u(.,t)=S^{A_1}(t-\gd)[u(.,\gd)]$. \smallskip

\noindent (ii) A function $u\in C((\gd,\infty);L^1(\BBR^N))$  is an extended semigroup solution of  $(\ref{A1})$ on $(\gd,\infty)$ if 
for any $t\geq\gt> \gd$, there holds $u(.,t)=S^{A_1}(t-\gt)[u(.,\gt)]$.

\es

\subsection{The Barenblatt-Prattle solutions}
We recall the explicit expression, due to Barenblatt and Prattle, of the solution $v=v_k$ of problem $(\ref{I-10})$.
If $p=2$
\bel{F-2}
v_k(x,t)=k(4\gp t)^{-\frac{N}{2}}e^{-\frac{|x|^2}{4t}},
\ee
and if $\frac{2N}{N+1}<p \ne 2$,
\bel{F-3}
v_k(x,t)=t^{-\gl}V\left(\frac{x}{t^{\frac{\gl}{N}}}\right),\text{ where }V(\xi)=\left(C_k-d|\xi|^{\frac{p}{p-1}}\right)_+^{\frac{p-1}{p-2}}
\ee
with
\bel{F-4}
\gl=\myfrac{N}{N(p-2)+p}\qquad\text {and }\quad d=\myfrac{p-2}{p}\left(\myfrac{\gl}{N}\right)^{\frac{1}{p-1}},
\ee
and where $C_k$ is connected to the mass $k$ by
\bel{F-5}
C_k=c(N,p)k^\ell \qquad\text {with }\quad\ell=\myfrac{p(p-2)\gl}{(p-1)N}.
\ee

The condition $p>\frac{2N}{N+1}$ appears in order $\gl$ be positive. Notice that, if $p>2$ then $d>0$, therefore the support of $v_k(.,t)$ is the ball 
$B_{\gd_k(t)}$ where $\gd_k(t)=\left(\frac{C_k}{d}\right)^{\frac{p-1}{p}}t^{\frac{\gl}{N}}$, while $v_k(x,t)>0$ for all $(x,t)\in Q_\infty$ if $\frac{2N}{N+1}<p<2$ (and also $p=2$ although the expression of $v_k$ is different). Furthermore, if $\frac{2N}{N+1}<p<2$, the limit of $v_k$ when $k\to\infty$ is explicit
\bel{F-5*}
v_\infty(x,t)=\Gl_N\left(\frac{t}{|x|^p}\right)^{\frac{1}{2-p}},
\ee
where $\Gl_N=(-d)^{\frac{p-1}{p-2}}$. This type of singular solution which is singular on the whole axis $(0,t)\subset Q_\infty$, is called a {\it razor blade} (see \cite{VaVe2} for some examples). To this solution corresponds a universal estimate.

\blemma{decay} Assume $1<p<2$ and let $v\in C(\overline{Q_\infty}\setminus B_{R_0}\ti \{0\})$ be a semigroup solution positive  of $(\ref{A1})$ 
\bel{Y1}\BA {l}
\prt_tv-\Gd_pv=0\qquad\text{in }Q_\infty
\EA\ee
which satisfies
\bel{Y2}\BA {l}
\Lim{t\to 0}\myint{K}{}v(x,t)dx=0,
\EA\ee
for any compact set $K\subset \BBR^N\setminus B_{R_0}$. Then there exists $c_1=c_1(N,p)>0$ such that
\bel{Y3}\BA {l}
\displaystyle\sup_{0\leq\gt\leq t}\myint{\{x:|x|>R\}}{}v(x,\gt)dx\leq c_1 \left(\myfrac{t}{(R-R_0)^{\frac{N}{\gl}}}\right)^{\frac{1}{2-p}}\qquad\forall R>R_0, t>0.
\EA\ee
If we assume moreover that ${\dsps \lim_{|x|\to\infty} v(t,x)=0}$ locally uniformly with respect to $t\geq 0$, then
\bel{Y3-1}\BA {l}
v(x,t)\leq \Gl_1\left(\myfrac{t}{(|x|-R_0)^p}\right)^{\frac{1}{2-p}}\qquad\forall (x,t)\in Q_\infty, |x|>R_0, 
\EA\ee
where $\Gl_1$ is the value of the constant in $(\ref{F-5*})$ when $N=1$.
\es
\Proof The first estimate is a consequence of
\bel{Y4}\BA {l}
\displaystyle\sup_{0\leq\gt\leq t}\myint{B_\gr(a)}{}v(x,t)dx\leq c_2 \left(\myint{B_{2\gr}(a)}{}v(x,0)dx+\left(\myfrac{t}{\gr^{\frac{N}{\gl}}}\right)^{\frac{1}{2-p}}\right)
\EA\ee
 in \cite[Lemma III.3.1]{DH} under the assumption that $v(.,0)$ is continuous with compact support. Actually this assumption is not used. In this proof the first step is the following estimate obtained by a suitable choice of test function:
\bel{Y5}\BA {l}
\displaystyle\sup_{0\leq\gt\leq t}\myint{B_R(a)}{}v(x,t)dx\leq \myint{B_{2R}(a)}{}v(x,0)dx+\myfrac{c_3}{R}\myint{0}{t}\myint{B_R(a)}{}|\nabla v|^{p-1}dx\,d\gt
\EA\ee
valid for any $a\in\BBR^N\setminus\{(0,0)\}$ and $R\leq |a|/2$. The second step to get $(\ref{Y4})$ is to estimate the integral on the right-hand side by relation (I.4.2) in \cite[Lemma I.4.1]{DH} with the same choice of $\ge$. We apply estimate $(\ref{Y4})$ with a sequence of points in a fixed direction ${\bf e}$ (with $|{\bf e}|=1$) $a=a_k=\left(2^k(R-R_0)+R_0\right){\bf e}$ and $\gr=\gr_k=2^{k-1}(R-R_0)$ (actually we start with $\gr<\gr_k$ and let it grow up to $\gr_k$). Then we get
\bel{Y6}\BA {l}
\displaystyle\sup_{0\leq\gt\leq t}\myint{B_{\gr_k}(a_k)}{}v(x,t)dx\leq c_4 2^{-\frac{N(k-1)}{\gl(2-p)}}\left(\myfrac{t}{(R-R_0)^{\frac{N}{\gl}}}\right)^{\frac{1}{2-p}}.
\EA\ee
Since the ball $B_{\gr_k}(a_k)$ and $B_{\gr_{k+1}}(a_{k+1})$ are overlapping there exist a finite number of points $\{{\bf e}_j\}^{d_1}_{j=1}$ and $\{{\bf e'}_j\}_{j=1}^{d_2}$ ($d_1$ and $d_2$ depend only on $N$) on the unit sphere such that 
$$\left\{x\in\BBR^N:\abs x\geq R\right\}\subset\left(\bigcup_{j=1}^{d_1}\bigcup_{k=1}^{\infty}B_{\gr_k}(2\gr_k{\bf e}_j)\right) \bigcup \left(\bigcup_{j=1}^{d_2}B_{\frac{R-R_0}{2}}(R{\bf e'}_j)\right). 
$$
Therefore
$$
\displaystyle\sup_{0\leq\gt\leq t}\myint{\{x:|x|>R\}}{}v(x,\gt)dx\leq c_4 \left[d_1\sum_{k=0}^{\infty}2^{-\frac{Nk}{\gl(2-p)}}+ d_2 2^\frac{N}{\gl(2-p)}\right]\left(\myfrac{t}{(R-R_0)^{\frac{N}{\gl}}}\right)^{\frac{1}{2-p}}
$$
which is $(\ref{Y3})$.\smallskip

Estimate $(\ref{Y3-1})$ follows from comparison with the 1-dim form of $v_\infty$
\bel{F-6-1}
v_\infty(s,t)=\Gl_1\left(\frac{t}{s^p}\right)^{\frac{1}{2-p}} \qquad\forall s,t>0.
\ee
For $\ge>0$, the function 
$$(x,t)\mapsto v_\infty(x_1-R_0-\ge,t)+\ge
$$
where $x=(x_1,...,x_N)=(x_1,x')$, is a solution of $(\ref{Y1})$ in $H_{1,R_0+\ge}\ti(0,\infty)$ where $H_{1,m}=\{x\in\BBR^N:x_1>m\}$. For $R$ large enough $v(x,t)\leq v_\infty(x_1+R_0+\ge,t)+\ge$ on the set $\left((H_{1,R_0+\ge}\cap \prt B_R)\cup (\prt H_{1,R_0+\ge}\cap B_R)\right)\ti [0,T]$ 
for any $T>0$, and for $t=0$. By the maximum principle $v(x,t)\leq v_\infty(x_1-R_0-\ge,t)+\ge$ in 
$(H_{1,R_0+\ge}\cap B_R)\ti (0,T]$. Letting successively $R\to\infty$, $T\to\infty$ and  $\ge\to 0$ and using the  invariance of the equation by rotation implies $(\ref{Y3-1})$.
\qeda\medskip

\bprop{stab} Let $p>\frac{2N}{N+1}$ and $\{v^{n}\}\subset C([0,\infty);L^1(\BBR^N))$ be a sequence of  positive semigroup solutions of $(\ref{Y1})$ on $(0,\infty)$ such that $v^{n}(.,0)$ has support in $B_{\ge_n}$ where $\ge_n\to 0$. If
$$
\myint{\BBR^N}{}v^{n}(x,0)dx=k_n\to k\quad\text{as }n\to\infty
$$
then $v^n\to v_k$ locally uniformly in $Q_\infty$.
\es
\Proof We first give the proof in the case $\frac{2N}{N+1}<p<2$. By {\it a priori}  estimates, up to a subsequence $v^n$ converges locally uniformly in $Q_\infty$ to a solution $v$ of $(\ref{Y1})$ in $Q_\infty$. By Herrero-Vazquez mass conservation property \cite[Theorem 2] {HeVa} (valid if $p>\frac{2N}{N+1}$)
$$\myint{\BBR^N}{}v^n(x,t)dx=\myint{\BBR^N}{}v^n(x,0)dx=k_n.
$$
By $(\ref{Y3-1})$
$$v^n(x,t)\leq \Gl_1\left(\frac{t}{(|x|-\ge_n)^p}\right)^{\frac{1}{2-p}}\qquad\forall t>0,\,\forall |x|>\ge_n.
$$
Since $\frac{p}{2-p}>N$, the function
$$x\mapsto \left(\frac{t}{(|x|-\ge_n)^p}\right)^{\frac{1}{2-p}}
$$
belongs to $L^1(\BBR^N\setminus B_{\gd})$, for any $\gd>\ge_n$. Since $v^n(x,t)\to v(x,t)$ uniformly in $B_{\gd}$, it follows by the dominated convergence theorem
\bel{Y3-3}
\lim_{n\to\infty}\myint{\BBR^N}{}v^n(x,t)dx=\myint{\BBR^N}{}v(x,t)dx=k.
\ee
Because $v$ is a positive solution with isolated singularity at $(0,0)$, it follows from \cite{CQW} that $v=v_k$, solution of $(\ref{I-10})$. \smallskip

When $p\geq 2$, the function $v_k(.,t)$ has a compact support $D_{k_n}(t)$ for any $t>0$ and $D_{k_n}(t)\subset B_{R_n(t)}$ where
\bel{Y3-4}
R_n(t)=\ge_n+c_5k_n^{\frac{p-2}{p}}t^{\frac{1}{N(p-2)+p}}\leq \ge^*+c_5k_*^{\frac{p-2}{p}}t^{\frac{1}{N(p-2)+p}}
\ee
 where $c_5=c_5(N,p)>0$, $\ge^*=\sup\{\ge_n;n\in\BBN\}$ and $k_*=\sup\{k_n;n\in\BBN\}$. Using Lebesgue dominating theorem we obtain again $(\ref{Y3-3})$.\qeda

\subsection{Fundamental solutions}

\medskip

The following lemma is fundamental.
\blemma{int} Assume  $p>\frac{2N}{N+1}$ and $f$ is a continuous nondecreasing function defined on $\BBR$ such that $f(0)=0$. Then, for any $k,R,T>0$, 
\bel{F-6}
\myint{1}{\infty}f(s)s^{-\frac{p(N+1)}{N}}ds<\infty\Longrightarrow f(v_k)\in L^1(B_R\ti (0,T)).
\ee
\es
\Proof The result is already proved in \cite{KV1} in the case $p>2$. It is probably known in the case $p=2$, but we have not found any reference. It appears to be new in the case $\frac{2N}{N+1}<p<2$. Without any loss of generality we can assume $R=T=1$.\smallskip
 
\noindent{\it Case 1}: $p=2$. By linearity we can assume that $k=(4\gp)^{\frac{N}{2}}$. Let
$$I=\int\int_{B_1\ti (0,1)}f(v_k) dx\,dt=\gw_N\myint{0}{1}\myint{0}{1}f\left(t^{-\frac{N}{2}}e^{-\frac{r^2}{4t}}\right)r^{N-1}dr\,dt.
$$
Set $s= t^{-\frac{N}{2}}e^{-\frac{r^2}{4t}}$, then
$$\BA {l}
I=2^{N-1}\gw_N\myint{0}{1}\myint{ t^{-\frac{N}{2}}e^{-\frac{1}{4t}}}{t^{-\frac{N}{2}}}
\left[-\ln s-\ln\left( t^{\frac{N}{2}}\right)\right]^{\frac{N-2}{2}}f(s)s^{-1}ds\,t^{\frac{N}{2}}dt\\[4mm]\phantom{I}
\leq 2^{N-1}\gw_N\myint{0}{1}\myint{e^{-\frac{1}{4t}}}{t^{-\frac{N}{2}}}
\left[-\ln s-\ln\left( t^{\frac{N}{2}}\right)\right]^{\frac{N-2}{2}}f(s)s^{-1}ds\,t^{\frac{N}{2}}dt\phantom{I}
\leq 2^{N-1}\gw_N (I_1+I_2)
\EA$$
where,
$$\BA {l}I_1=\myint{0}{e^{-\frac{1}{4}}}\myint{0}{-\frac{1}{4\ln s}}\left[-\ln s-\ln\left( t^{\frac{N}{2}}\right)\right]^{\frac{N-2}{2}}\,t^{\frac{N}{2}}dt\,s^{-1}f(s)ds\\[4mm]\phantom{I_1}
=\myfrac{2}{N}\myint{0}{e^{-\frac{1}{4}}}\myint{0}{\frac{s}{(-4\ln s)^{\frac{N}{2}}}}(-\ln\gt)^{\frac{N-2}{2}}
\gt^{\frac{2}{N}}d\gt\,s^{-2-\frac{2}{N}}f(s)ds,\EA$$
by setting $\gt=st^{\frac{N}{2}}$. But
$$\BA {l}\myint{0}{\frac{s}{(-4\ln s)^{\frac{N}{2}}}}(-\ln\gt)^{\frac{N-2}{2}}
\gt^{\frac{2}{N}}d\gt\leq c_6\left[(-\ln\gt)^{\frac{N-2}{2}}
\gt^{\frac{N+2}{N}}\right]_0^{\frac{s}{(-4\ln s)^{\frac{N}{2}}}}\\[4mm]
\phantom{\myint{0}{\frac{s}{(-4\ln s)^{\frac{N}{2}}}}(-\ln\gt)^{\frac{N-2}{2}}
\gt^{\frac{2}{N}}d\gt}
\leq c_6s^{1+\frac{2}{N}}(-\ln s)^{-2}\left(1+\frac{N}{2}\frac{\ln(-4\ln s)}{-\ln s}\right)^{\frac{N-2}{2}}
\\[4mm]
\phantom{\myint{0}{\frac{s}{(-4\ln s)^{\frac{N}{2}}}}(-\ln\gt)^{\frac{N-2}{2}}
\gt^{\frac{2}{N}}d\gt}
\leq c_7s^{1+\frac{2}{N}}(-\ln s)^{-2},
\EA$$
thus 
$$I_1\leq c_8\myint{0}{e^{-\frac{1}{4}}}s^{-1}(-\ln s)^{-2}f(s)ds<\infty
$$
by Duhamel's rule. Further
$$\BA {l}I_2\leq\myint{e^{-\frac{1}{4}}}{\infty}\myint{0}{s^{-\frac{2}{N}}}\left[-\ln s-\ln\left( t^{\frac{N}{2}}\right)\right]^{\frac{N-2}{2}}\,t^{\frac{N}{2}}dt\,s^{-1}f(s)ds\\[4mm]\phantom{I_2}
\leq \myfrac{2}{N}\myint{e^{-\frac{1}{4}}}{\infty}\myint{0}{1}(-\ln\gt)^{\frac{N-2}{2}}
\gt^{\frac{2}{N}}d\gt\,s^{-2-\frac{2}{N}}f(s)ds\\[4mm]\phantom{I_2}
\leq c_9\myint{e^{-\frac{1}{4}}}{\infty}s^{-2-\frac{2}{N}}f(s)ds,
\EA$$
for some $c_9=c_9(N)>0$.
This implies the claim when $p=2$.\smallskip

\noindent{\it Case 2}: $\frac{2N}{N+1}<p<2$. We set $d^*=-d$. By rescaling we can assume that $C_k=d^*=1$. Therefore
$$I=\int\int_{B_1\ti (0,1)}f(v_k) dxdt=\gw_N\myint{0}{1}\myint{0}{1}f\left(t^{-\gl}\left[1+\left(\frac{r}{t^{\frac{\gl}{N}}}\right)^{\frac{p}{p-1}}\right]^{\frac{p-1}{p-2}}\right)r^{N-1}dr\,dt.
$$
Set $s=t^{-\gl}\left[1+\left(\frac{r}{t^{\frac{\gl}{N}}}\right)^{\frac{p}{p-1}}\right]^{\frac{p-1}{p-2}}$, then
$r=t^{\frac{\gl}{N}}\left[(t^\gl s)^{\frac{p-2}{p-1}}-1\right]^{\frac{p-1}{p}}$ and
$$\BA {l}
I=\frac{2-p}{p}\gw_N\myint{0}{1}\myint{t^{-\gl}(1+t^{-\frac{\gl p}{p-1}})^{\frac{p-1}{p-2}}}{t^{-\gl}}
(t^\gl s)^{-\frac{1}{p-1}}\left(\left(t^\gl s\right)^{\frac{p-2}{p-1}}-1\right)^{\frac{N(p-1)}{p}-1}f(s)
ds\,t^{2\gl} dt\\[3mm]\phantom{I}
=\frac{2-p}{p}\gw_N(I_1+I_2)\EA$$
where
$$\BA {l}
I_1=\myint{2^{\frac{p-1}{p-2}}}{1}\myint{a(s)}{1}(t^\gl s)^{-\frac{1}{p-1}}\left(\left(t^\gl s\right)^{\frac{p-2}{p-1}}-1\right)^{\frac{N(p-1)}{p}-1}t^{2\gl} dt\,f(s)
ds
\EA$$
$$\BA {l}
I_2=\myint{1}{\infty}\myint{a(s)}{s^{-\frac{1}{\gl}}}(t^\gl s)^{-\frac{1}{p-1}}\left(\left(t^\gl s\right)^{\frac{p-2}{p-1}}-1\right)^{\frac{N(p-1)}{p}-1}t^{2\gl} dt\,f(s)
\EA$$
and $a(s)$ is the inverse function of  $t\mapsto t^{-\gl}(1+t^{-\frac{\gl p}{p-1}})^{\frac{p-1}{p-2}}$. Clearly
$$ t^{-\gl}(1+t^{-\frac{\gl p}{p-1}})^{\frac{p-1}{p-2}}\leq t^{\frac{2\gl(p-1)}{2-p}}
\Longrightarrow a(s)\geq s^{\frac{2-p}{2\gl(p-1)}}.$$
Therefore
$$\BA {l}
I_1\leq \myint{2^{\frac{p-1}{p-2}}}{1}\myint{s^{\frac{2-p}{2\gl(p-1)}}}{1}(t^\gl s)^{-\frac{1}{p-1}}\left(\left(t^\gl s\right)^{\frac{p-2}{p-1}}-1\right)^{\frac{N(p-1)}{p}-1}t^{2\gl} dt\,f(s)ds\\[4mm]\phantom{I_1}
\leq \myfrac{1}{\gl}\myint{2^{\frac{p-1}{p-2}}}{1}
\myint{s^{\frac{p}{2(p-1)}}}{s}(1-\gt^{\frac{2-p}{p-1}})^{\frac{N(p-1)}{p}-1}\gt^{\frac{1}{\gl}+\frac{N(p-2)}{p}}d\gt\,
s^{-2-\frac{1}{\gl}}f(s)ds.
\EA$$
Since $\frac{1}{\gl}+\frac{N(p-2)}{p}>-1$ and $\frac{N(p-1)}{p}-1>-1$,
$$\myint{s^{\frac{p}{2(p-1)}}}{s}(1-\gt^{\frac{2-p}{p-1}})^{\frac{N(p-1)}{p}-1}\gt^{\frac{1}{\gl}+\frac{N(p-2)}{p}}d\gt
<\myint{0}{1}(1-\gt^{\frac{2-p}{p-1}})^{\frac{N(p-1)}{p}-1}\gt^{\frac{1}{\gl}+\frac{N(p-2)}{p}}d\gt<\infty.
$$
Furthermore $-2-\frac{1}{\gl}=-p-\frac{p}{N}$ thus
$$I_1\leq c_{10}\myint{2^{\frac{p-1}{p-2}}}{1}f(s)s^{-\frac{p(N+1)}{N}}ds.
$$
We perform the same change of variable with $I_2$
$$\BA {l}I_2\leq \myint{1}{\infty}\myint{s^{\frac{2-p}{2\gl(p-1)}}}{s^{-\frac{1}{\gl}}}
(t^\gl s)^{-\frac{1}{p-1}}\left(\left(t^\gl s\right)^{\frac{p-2}{p-1}}-1\right)^{\frac{N(p-1)}{p}-1}t^{2\gl} dt\,f(s)ds\\[4mm]\phantom{I_2}
\leq\myfrac{1}{\gl}\myint{1}{\infty}\myint{s^{\frac{p}{2(p-1)}}}{1}(1-\gt^{\frac{2-p}{p-1}})^{\frac{N(p-1)}{p}-1}\gt^{\frac{1}{\gl}+\frac{N(p-2)}{p}}d\gt\,
s^{-2-\frac{1}{\gl}}f(s)ds.
\EA$$
Again
$$\myint{s^{\frac{p}{2(p-1)}}}{1}(1-\gt^{\frac{2-p}{p-1}})^{\frac{N(p-1)}{p}-1}\gt^{1+\frac{N(p-2)}{p}}d\gt
<\myint{0}{1}(1-\gt^{\frac{2-p}{p-1}})^{\frac{N(p-1)}{p}-1}\gt^{\frac{1}{\gl}+\frac{N(p-2)}{p}}d\gt<\infty,$$ 
then
$$I_2\leq c_{11}\myint{1}{\infty}f(s)s^{-\frac{p(N+1)}{N}}ds.
$$
Therefore $(\ref{F-6})$ holds.\qeda\medskip

Notice that the assumption implies that $v_k\in C(Q_\infty)\cap L^\infty(\gd,\infty;L^1(\BBR^N)\cap L^\infty(\BBR^N))$ for every $\gd>0$.\medskip

\noindent{\bf Proof of \rth{Fund}.} {\it Existence.} Let $\ge>0$, $Q_{\ge,\infty}=\BBR^N\ti(\ge,\infty)$ and denote by $u_\ge$ the solution of 
\bel{B1}\left\{\BA {ll}
\!\prt_tu-\Gd_pu+f(u)=0\qquad&\text{in }Q_{\ge,\infty}\\\phantom{\Gd_pu+f(u)}
u(.,\ge)=v_k(.,\ge)\qquad&\text{in }\BBR^N.
\EA\right.\ee
Since $v_k(.,\ge)$ is a smooth positive function belonging to $L^1(\BBR^N)$ the function $u_\ge$ is constructed by truncation. By the maximum principle
\bel{B2}
u_\ge(x,t+\ge)\leq v_k(x,t+\ge)\qquad\forall (x,t)\in Q_{\ge,\infty}.
\ee
For $0<\ge'<\ge$, $u_{\ge'}(x,\ge)\leq v_k(x,\ge)=u_\ge(x,\ge)$, thus  $u_{\ge'}(x,t+\ge)\leq u_\ge(x,t+\ge)$ in $Q_{\ge,\infty}$. Set $\tilde u=\lim_{\ge\to 0}u_\ge$, then $\tilde u\leq v_k$ in $Q_\infty$. By the standard local regularity theory for degenerate equations, $\nabla u_\ge$ remains locally compact in  $(C^1_{loc}(Q_\infty))^N$, thus $\tilde u$ satisfies $(\ref{A1})$ in $Q_\infty$.  \smallskip

In order to prove that
$$\frac{d}{dt}\myint{\BBR^N}{}u_{\ge}(x,s)dx
+\myint{\BBR^N}{}f(u_{\ge}(x,s))dx=0
$$
we recall that $u_\ge$ can be obtained as the limit of thru the iterative implicit scheme $(\ref{X4})$ with $q\in [1,\infty]$ is arbitrary since $u_{\ge,0}\in L^1(\BBR^N)\cap L^\infty(\BBR^N)$.  For $h>0$ we can write it under the form
$$u_{\ge,i}-h\Gd_pu_{\ge,i}=-hf(u_{\ge,i})+u_{\ge,i-1}.
$$
By $(\ref{X6})$, and denoting by $\tilde U_{\ge,h}$ the piecewise constant function such that 
$\tilde U_{\ge,h}(jh)=u_{\ge,j}$, we obtain since $u_{\ge,0}=v_k(\ge)$
\bel{B4}
\myint{\BBR^N}{}(u_{\ge,i}-v_k(\ge))(x)dx=-\myint{\ge}{ih}\myint{\BBR^N}{}f(\tilde U_{\ge,h}(x))dx dt.
\ee
Letting $h\to 0$ and $ i\to\infty$ such that $ih=t>\ge$ and using the uniform convergence, we obtain
\bel{B5}
\myint{\BBR^N}{}u_\ge(x,t)dx-\myint{\BBR^N}{}v_k(x,\ge)dx=-\myint{\ge}{t}\myint{\BBR^N}{}f(u_\ge(x,s))dx dt.
\ee
Since $0\leq u_\ge\leq v_k$ and $v_k(.,t)$ has constant mass equal to $k$, we derive
\bel{B6}
\abs{\myint{\BBR^N}{}u_\ge(x,t)dx-k}\leq \myint{\ge}{t}\myint{\BBR^N}{}f(v_k(x,s))dx dt.
\ee
Because $f(v_k) \in L^1(\BBR^N\ti (0,T))$, we can let $\ge\to 0$, using the monotone convergence theorem, in order to get 
\bel{B6*}
\abs{\myint{\BBR^N}{}u(x,t)dx-k}\leq \myint{0}{t}\myint{\BBR^N}{}f(v_k(x,s))dx dt.
\ee
This implies that
\bel{B7}
\lim_{t\to 0}\myint{\BBR^N}{}u(x,t)dx=k.
\ee
If $\gf\in C_c(\BBR^N)$, let $\gz\in C_c^\infty(\BBR^N) $ such that $0\leq\gz\leq 1$, $\gz=1$ on the support of $\gf$ and $\gz(0)=1$. Then 
$$\BA{l}\myint{\BBR^N}{}u(x,t)\gf(x)dx=\myint{\BBR^N}{}u(x,t)\gf(x)\gz(x)dx\\[4mm]\phantom{\myint{\BBR^N}{}u(x,t)\gf(x)dx}=
\gf(0)\myint{\BBR^N}{}u(x,t)dx+\myint{\BBR^N}{}u(x,t)(\gf(x)\gz(x)-\gf(0))dx.
\EA$$
Thus
$$\abs{\myint{\BBR^N}{}u(x,t)\gf(x)dx-\gf(0)\myint{\BBR^N}{}u(x,t)dx}\leq \myint{\BBR^N}{}v_k(x,t)\abs{\gf(x)\gz(x)-\gf(0)}dx.
$$
Because $\abs{\gf(x)\gz(x)-\gf(0)}$ is continuous and vanishes at zero and $v_k(.,0)=k\gd_0$, it follows from $(\ref{B7})$
\bel{B8}
\lim_{t\to 0}\myint{\BBR^N}{}u(x,t)\gf(x)dx=k\gf(0).
\ee

\noindent{\it Uniqueness}. The proof uses some ideas from \cite[Th 2.4]{KV1}. Assume $\tilde u$ is any nonnegative solution of problem $(\ref{I-9})$, then, for any $\ge>0$ we denote by $\tilde v_\ge$ the solution of 
\bel{B9}\left\{\BA {ll}
\prt_tv-\Gd_pv=0\qquad&\text{in } Q_{\ge,\infty}\\\phantom{,--}
v(.,\ge)=\tilde u(.,\ge)\qquad&\text{in } \BBR^N.
\EA\right.\ee
By the maximum principle $\tilde v_\ge\geq \tilde u$ in $Q_{\ge,\infty}$. When $\ge\to 0$, $\tilde v_\ge$ converges locally uniformly to a solution $\tilde v$ of the same equation in  $Q_{\infty}$. Furthermore, using again \cite[Lemma 2]{HeVa},
$$\myint{\BBR^N}{}\tilde v_\ge(x,t+\ge)dx=\myint{\BBR^N}{}\tilde u(x,\ge)dx.$$
By Fatou's Lemma and using the fact that 
$$\lim_{\ge\to 0}\myint{\BBR^N}{}\tilde u(x,\ge)dx=k,
$$
we derive
\bel{B10}
\myint{\BBR^N}{}\tilde v(x,t)dx\leq k.
\ee
Since $\tilde v\geq \tilde u$, equality holds in $(\ref{B10})$. Since the fundamental solution is unique \cite[Th 4.1]{CQW}, it implies $\tilde v=v_k$ and $\tilde u\leq v_k$. We end the proof as in \cite[Th 4.1]{CQW}, using the $L^1$-contraction mapping principle and the fact that any solution of $(\ref{I-9})$ is smaller than $v_k$: for $t>s>0$, there holds
\bel{B11}\BA {l}
\myint{\BBR^N}{}\abs{u(x,t)-\tilde u(x,t)}dx\leq \myint{\BBR^N}{}\abs{u(x,s)-\tilde u(x,s)}dx\\[4mm]
\phantom{\myint{\BBR^N}{}\abs{u(x,t)-\tilde u(x,t)}dx}
\leq \myint{\BBR^N}{}\abs{u(x,s)-v_k(x,s)}dx+\myint{\BBR^N}{}\abs{v_k(x,s)-\tilde u(x,s)}dx\\[4mm]
\phantom{\myint{\BBR^N}{}\abs{u(x,t)-\tilde u(x,t)}dx}
\leq \myint{\BBR^N}{}(v_k(x,s)-u(x,s))dx+\myint{\BBR^N}{}(v_k(x,s)-\tilde u(x,s))dx.
\EA\ee
When $s\to 0$ the right-hand side of the last line goes to $0$. This implies the claim.\qeda
\medskip

The next result shows some geometric properties of the $u_k$.

\bprop{Rad} The solution $u=u_k$ of problem $(\ref{I-6})$ is radial and nonincreasing with respect to $|x|$.
\es
\Proof It is sufficient to prove the result with the approximation $u_\ge(.,t)$. By $(\ref{F-3})$, $v_k(.,t)$ is radial and decreasing, therefore $u_\ge(.,t)$ is radial too by uniqueness. We notice that $u_\ge$ is the increasing limit, when $R\to\infty$, of the solution $u_{\ge,R}$ of 
\bel{B3}\left\{\BA {ll}
\!\prt_tu-\Gd_pu+f(u)=0\qquad&\text{in }Q^{B_R}_{\ge,\infty}\\\phantom{\Gd_pu+f(u)--}
u=0\qquad&\text{in }\prt B_R\ti(\ge,\infty)\\\phantom{\Gd_pu+f(u)}
u(.,\ge)=v_k(.,\ge)\qquad&\text{in }B_R.
\EA\right.\ee
For $\gl\in (0,R)$, we set $\Gs_\gl=B_R\cap \{x=(2\gl-x_1,x'):x_1>\gl\}\cap B_R$ and define $w_\gl$ by
$$w_\gl(x,t)=w_\gl(x_1,x',t):=u_{\gl,\ge,R}(x)-u_{\ge,R}(x)=u_{\ge,R}(2\gl-x_1,x',t)-u_{\ge,R}(x_1,x',t).$$
If $Q^{\Gs_\gl}_{\ge,\infty}=\Gs_\gl\ti (\ge,\infty)$, there holds
\bel{B3*}\left\{\BA {ll}
\!\prt_tw_\gl+\CA w_\gl+d(x)w_\gl=0\qquad&\text{in }Q^{\Gs_\gl}_{\ge,\infty}\\\phantom{\Gd_pw_\gl+f(w_\gl)--}
w_\gl\geq 0\qquad&\text{in }\prt \Gs_\gl\ti(\ge,\infty)\\\phantom{\Gd_pw_\gl+f(w_\gl)}
w_\gl(.,\ge)\geq 0\qquad&\text{in }\Gs_\gl.
\EA\right.\ee
where
$$d(x)=\left\{\BA {ll}\frac{f(u_{\gl,\ge,R})-f(u_{\ge,R})}{u_{\gl,\ge,R}-u_{\ge,R}}\quad&\text{if }u_{\gl,\ge,R}\neq u_{\ge,R}\\
0\quad&\text{if }u_{\gl,\ge,R}= u_{\ge,R}
\EA\right.$$
and
$$\CA w_\gl=-\Gd_pu_{\gl,\ge,R}+\Gd_pu_{\ge,R}.
$$
Notice that $d\geq 0$ since $f$ is nondecreasing and $\CA$ is elliptic \cite[Lemma 1.3]{FV}. Furthermore the boundary data are continuous, therefore $w_\gl\geq 0$. Letting $\gl\to 0$, changing $\gl$ in $-\gl$ and replacing the $x_1$ direction, by any direction going thru $0$, we derive that 
$u_{\ge,R}(.,t)$ is radially decreasing. Letting $R\to\infty$ yields to $u_{\ge}(.,t)$ is radially decreasing too. \qeda\medskip

In the next result we characterize positive solutions of $(\ref{A1})$ with an isolated singularity at $t=0$

\bprop{isolsing} Assume $p>\frac{2N}{N+1}$ and $f$ is  continuous nondecreasing function vanishing only at $0$ and  satisfying $(\ref{CFS})$. If $u\in C(\overline{Q_\infty}\setminus \{(0,0)\})$ is a positive semigroup solution of $(\ref{A1})$ in $Q_\infty$ such that $u(x,0)=0$, for all $x\neq 0$ and
	$$ \lim_{t \to 0}\myint{\BBR^N}{}u(x,t)dx < \ity, $$
then there exists $k\geq 0$ such that $u=u_k$.
\es
\Proof  Using \cite[Lemma 2.2 ]{KV2} when $p\geq 2$, or the proof of \rth{Fund} when $\frac{2N}{N+1}<p<2$ jointly with the fact that
$$t\mapsto\myint{\BBR^N}{}u(x,t)dx
$$
is decreasing, we derive that $u\leq v_m$ for some $m\geq 0$ and there exists $k\geq 0$ such that 
$$\lim_{t\to 0}\myint{\BBR^N}{}u(x,t)dx=k.
$$
Since $u(.,0)$ vanishes if $x\neq 0$, it implies
$$\lim_{t\to 0}\myint{\BBR^N}{}u(x,t)\gf(x)dx=k\gf(0)\qquad\forall\gf\in C_c(\BBR^N).
$$
Therefore $u$ satisfies $(\ref{I-9})$. By uniqueness, $u=u_k$.\qeda \medskip

\subsection{Strong singularities}
This section is devoted to study the limit of the sequence of the solutions $u_k$ to $(\ref{I-9})$ as $k \to \ity$ with $f(s)=s^\ga\ln^\gb(s+1)$ where $p>2$, $\ga \in [1,p-1)$ and $\gb>0$. \medskip

\noindent{\bf Proof of \rth{IS1}.} By the comparison principle, 
$$ u_k(x,t) \leq v_k(x,t) \leq c_{12}k^\frac{(p-1)\ell}{p-2}t^{-\gl}, $$ 
where $v_k$ is the solution of $(\ref{I-10})$ in $Q_\ity$ and $c_{12}=c_{12}(N,p)>0$ in $(\ref{F-5})$. We set 
	\bel{gth}\gth_k(t)=c_{12}^{\ga-1}k^{\frac{\ell(\ga-1)(p-1)}{p-2}}t^{-\gl (\ga-1)}\ln^\gb(c_{12}k^{\frac{(p-1)\ell}{p-2}}t^{-\gl}+1)\ee 
then 
	\bel{ine1} \prt_t u_k - \Gd_p u_k + u_k \gth_k(t) \geq 0. \ee
Next we write $u_k(x,t)=b_k(t)w_k(x,s_k(t))$ (the functions $b_k$ and $s_k$ will be defined later). For simplicity, we drop the subscript $k$ in $b_k$ and $s_k$. Inserting in $(\ref{ine1})$, we get
	\bel{ine2} b^{2-p}(t)s\pr (t) \prt_{s} w_k(x,s) - \Gd_p w_k(x,s) + b^{1-p}[b\pr (t)+b(t)\gth_k(t)]w_k(x,s) \geq 0. \ee
We choose the functions $b$ and $s$ such that
	$$ b^{2-p}(t)s\pr (t)=1 \qq \text{and} \qq b\pr(t) + b(t)\gth_k(t)=0, $$
which implies 
	\bel{b,s} b(t)=\exp\big(-\myint{0}{t}\gth_k(\tau)d\tau\big) \q \text{ and } \q
	s(t)=\myint{0}{t}\exp\big(-(p-2)\myint{0}{\tau}\gth_k(\sigma)d\sigma\big)d\tau. \ee
Then $\prt_{s} w_k - \Gd_p w_k \geq 0$ in $\BBR^N \ti (0,s_{k,0})$ with some $s_{k,0}>0$ and $w_k(.,0)=k\gd_0$. It follows by comparison principle that $w_k \geq v_k$ in $\BBR^N \ti (0,s_{k,0})$.
Hence 
	\bel{ine3} u_k(x,t) \geq b(t)v_k(x,s)=\exp\left(-\myint{0}{t}\gth_k(\tau)d\tau\right)s^{-\gl}\big(c_{13}k^\ell-c_{14}s^\frac{-p\gl}{(p-1)N}\abs{x}^\frac{p}{p-1}\big)_+^\frac{p-1}{p-2}. \ee
Let $\gd_1> \frac{\ell (\ga-1)(p-1)}{p-2}$ and $0<\gd_2<1-\gl (\ga-1)$. Using $(\ref{gth})$ there exists $t_0>0$ depending on $\gd_1$, $\gd_2$ and $k$ large enough, such that, for any $t \in (0,t_0)$ there holds
	\bel{est-gth} \myint{0}{t}\gth_k(\tau)d\tau \leq c_{15}  k^{\gd_1}t^{\gd_2} \forevery t \in (0,t_0) \ee
with $c_{15}=c_{15}(c_i,\ga,\gb,p,N)>0$. It follows from $(\ref{b,s})$ and $(\ref{est-gth})$ that
	\bel{est-s} t\,\exp\bigg[-(p-2)c_{15}  k^{\gd_1}t^{\gd_2}\bigg] \leq s(t) \leq t. 
	\ee
Since $J<\infty$ holds, there exists the solution $\phi_\ity$ of $(\ref{I-1})$. The sequence $\{u_k\}$ is increasing and is bounded from above by $\phi_\ity$, then the function $\unl U(x,t):={\dsps \lim_{k \to \ity}u_k(x,t)}$ satisfies $\unl U(x,t) \leq \phi_\ity(t)$ for every $(x,t) \in Q_\ity$. We restrict $x \in B_1$ and we choose $t$ such that
	\bel{est1-k} c_{13}k^\ell - c_{14}s(t)^\frac{-p\gl}{(p-1)N}>\myfrac{1}{2}c_{13}k^\ell \Llra k > \bigg(\myfrac{2c_{14}}{c_{13}}\bigg)^\frac{1}{\ell}s(t)^{-\frac{1}{p-2}}. \ee
By $(\ref{est-s})$, we only need to choose $t$ such that
	\bel{est2-k} k \geq \bigg(\myfrac{2c_{14}}{c_{13}}\bigg)^\frac{1}{\ell} t^\frac{-1}{p-2}\,\exp\bigg(c_{15}  k^{\gd_1}t^{\gd_2}\bigg). \ee 
We choose $t$ under the form 
	\bel{k} t=k^{-\frac{1}{\gg}} \text{ with } \gg>0, \ee
then $(\ref{est2-k})$ becomes
	\bel{est1-t} t^{-\gg} \geq \bigg(\myfrac{2c_{14}}{c_{13}}\bigg)^\frac{1}{\ell} t^\frac{-1}{p-2} \exp\left( c_{15}  t^{\gd_2-\gd_1 \gg} \right). \ee
In order to obtain $(\ref{est1-t})$, it is sufficient to choose $\gg$ such that
	\bel{con-gamma} \myfrac{1}{p-2} < \gg < \myfrac{\gd_2}{\gd_1}. \ee
Indeed, since $\ga<p-1$, we may choose $\gd_1$ and $\gd_2$ close enough $\frac{\ell (\ga-1)(p-1)}{p-2}$ and $1-\gl (\ga-1)$ respectively such that $(\ref{con-gamma})$ holds true. When $t$ has the form $(\ref{k})$ where $\gg$ satisfies $(\ref{con-gamma})$, from $(\ref{ine3})$, $(\ref{est-gth})$-$(\ref{est1-k})$ and the fact that $\unl U \geq u_k$ in $Q_\ity$, we deduce that
	\bel{est1-U} \unl U(x,t) \geq c_{16}t^{-\gl}\exp\bigg[c_{17}\ln (t^{-1})- c_{15}t^{\gd_2-\gd_1 \gg}\bigg]    
	\ee
for every $(x,t) \in B_1 \ti (0,t_0)$ with $t_0$ small enough and $c_{16}=c_{16}(N,p)$, $c_{17}=c_{17}(N,p,\gg)$. Since $\gg$ satisfies $(\ref{con-gamma})$, 
	$$ c_{17}\ln (t^{-1})- c_{15}t^{\gd_2-\gd_1 \gg} >0 $$
for every $t \in (0,t_0)$.
Therefore $ \lim_{t \to 0}\unl U(x,t) = \ity $ uniformly with respect to $x \in B_1$. We next proceed as in \cite[Lemma 3.1]{VaVe2} to deduce that $\unl U(x,t)$ is independent of $x$ and therefore it is a solution of $(\ref{I-1})$.
Since $J<\ity$, $\unl U(x,t) = \gf_\ity(t)$ for every $(x,t) \in Q_\ity$. \qeda \medskip

\rth{IS2} is proved by the same arguments as \rth{IS1}, using the fact that $\unl U(x,t)$ is independent of $x$.

\mysection{Non-Uniqueness}
The next result shows that $K=\infty$ is the necessary and sufficent condition so that a local solution of 
	\bel{Con1} (r^{N-1}{\abs {w\pr}}^{p-2}w\pr)\pr = r^{N-1}f(w) \ee
can be continued as a global solution. More precisely,
\blemma{Con} Every positive and increasing solution of $(\ref{Con1})$ defined in an interval $[a,a^*]$ to the right of $a > 0$ can be continued as a solution of $(\ref{Con1})$ on $[a,+\ity)$ if and only if $f$ satisfies 
	\bel{K-inf-2} \myint{\ga}{\ity}\frac{ds}{(F(s))^{1/p}}=\ity \ee
for any $\ga >0$.
\es
\Proof The proof is is an extension to the case $p\neq 2$ of the one of \cite[Lemma 2.1]{VaVe1} for the case $p=2$. \smallskip

\noindent{\it Step 1.} We first assume that $w$ is defined on a maximal interval $[a,a^*)$ with $a^*<\ity$ and ${\dsps \lim_{r \to a^*}w(r)=+\ity}$. Since $w$ is a nondecreasing function, $w' \geq 0$. And hence we may write $(\ref{Con1})$ under the following form
	$$ \frac{N-1}{r}(w')^{p-1}+(p-1)(w')^{p-2}w\ppr = f(w), $$
which implies that
	\bel{Conti3} (p-1)(w')^{p-2}w\ppr \leq f(w) \ee
and hence 
	$$ \frac{p-1}{p}({w\pr}^p)\pr \leq (F(w))\pr. $$
Taking the integral over $[a,r]$, we get
	$$ \frac{p-1}{p}[(w')^p(r) - (w')^p(a)] \leq F(w(r)) - F(w(a)) \leq F(w(r)). $$
Since $f$ is positive on $(0,\infty)$, $F(s)\to\infty$ when $s\to\infty$, thus there exists $\tilde a\in (a,a^*)$ such that 
	$$ 0<w\pr(r))^p \leq \frac{p}{2(p-1)}F(w(r))\Longrightarrow \frac{w\pr(r)}{F(w(r))^{1/p}} \leq  \Big(\frac{p}{2(p-1)}\Big)^{1/p}\qquad\forall r\in [\tilde a,a^*). $$
Taking the integral over $ [\tilde a,r]$, we obtain
	$$ \myint{w(\tilde a)}{w(r)}\frac{ds}{F(s)^{1/p}} \leq  \Big(\frac{p}{2(p-1)}\Big)^{1/p}(r-\tilde a). $$
Letting $r \to a^*$ yields to
	$$ \myint{w(\tilde a)}{\ity}\frac{ds}{F(s)^{1/p}} \leq \Big(\frac{p}{2(p-1)}\Big)^{1/p}(a^*-\tilde a) < \ity $$
and $(\ref{K-inf-2})$ is not satisfied. \smallskip

\noindent{\it Step 2.} We assume that
	$$ \myint{\ga}{\ity}\frac{ds}{F(s)^{1/p}}  < \ity $$
for some $\ga >0$, and we fix $A>a$. By \cite[Theorem 1]{Va} there exists a function $\gg$ defined on $(a,A)$ such that
	$$ w(r) < \gg(r) \forevery r \in (a,A) $$
for any solution of $(\ref{Con1})$ on $(a,A)$. Moreover, $\gg$ can be assumed convex, and ${\dsps \lim_{t \to a}}\gg(r)={\dsps\lim_{r \to A}\gg(r)}=+\ity$. If $w$ is a solution of $(\ref{Con1})$ on $(a,a+\ge)$ such that $w(a) > {\dsps \min_{a < r <A}\gg(r)}$ and $\gg\pr(a)>0$, it is clear that $w(r^*)=\gg(r^*)$ for some $r^*<A$ and $w(r) > \gg(r)$ for $r \in (r^*,r^*+\ge)$, so $w$ can not be defined on the whole $(a,A)$, and there exists $a^*<A$ such that ${\dsps \lim_{r \to a^*}w(r)=\ity}$. \qeda \medskip

\noindent{\bf Proof of \rth{glob}} By the Picard-Lipschitz fixed point theorem  in the case $1<p<2$ and \cite[Th 5.2]{GV} in the case $p\geq 2$, there exists a unique solution $w_a$ to $(\ref{I-7})$ defined on a maximal interval $[0,r_a)$ and  $w_a$ is an increasing function. Since Keller-Osserman estimate does not hold, by \rlemma{Con}, the solution can be continued on the whole $[0,+\ity)$ and global uniqueness follows from the local uniqueness. The function $r \mt w_a(r)$ is increasing and
	$$ w_a(r) \geq a+\frac{p-1}{p}\Big(\frac{f(a)}{N}\Big)^{\frac{1}{p-1}}r^\frac{p}{p-1} \q \text{and} \q 
	 w_a\pr (r)\geq  \Big(\frac{f(a)}{N}\Big)^{\frac{1}{p-1}}r^\frac{1}{p-1} $$
for any $r > 0$. \qeda 

\bprop{Ex2} Assume $p>2N/(N+1)$, $f$ is locally Lipschitz continuous and $K=\infty$ hold. For any positive function $u_0 \in C(Q_\ity)$ which satisfies
	\be w_a(\abs x)\leq u_0(x) \leq w_b(\abs x) \q \forall x \in {\BBR}^N	\ee
for some $0<a<b$, there exists a positive function $\ovl u \in C(\ovl {Q_\ity})$ solution of $(\ref{A1})$ in $Q_\ity$ and satisfying $\ovl u(.,0)=u_0$ in $\BBR^N$. Furthermore
	\bel{Ex2.1} w_a(\abs x)\leq \ovl u(x,t) \leq w_b(\abs x) \q \forall (x,t)\in Q_\ity.	\ee
\es
\Proof Clearly $w_a$ and $w_b$ are ordered solutions of $(\ref{A1})$. We denote by $u_n$ the solution to the initial-boundary problem
	\be \left\{ \BA{lll}	
	\prt_t u_n -\Gd_p u_n + f(u_n) &= 0  \q &\text{in } Q_n: = B_n \ti (0,\ity) \\ \phantom{------,}
	u_n(x,0)&=u_0(x) &\text{in } B_n \\ \phantom{------,}
	u_n(x,t)&=(w_a(\abs x)+w_b(\abs x))/2  &\text{in } \prt B_n \ti (0,\ity).
	\EA \right. \ee
By the maximum principle, $u_n$ satisfies $(\ref{Ex2.1})$ in $Q_n$. Using locally parabolic equation regularity \cite[Th 1.1, chap III]{Di} if $p \geq 2$ or \cite[Th 1.1, chap IV]{Di} if $1<p<2$, we derive that the set of functions $\{u_n\}$ is eventually equicontinuous on any compact subset of $\ovl {Q_\ity}$. Using a diagonal sequence, combined with \rprop{reg-Mar}, we conclude that there exists a subsequence $\{u_{n_k}\}$ which converges locally uniformly in $\ovl {Q_\ity}$ to some weak solution $\ovl u \in C(\ovl {Q_\ity})$ which has the desired properties. \qeda
\bprop{Ex3} Assume $p>2N/(N+1)$, $f$ is locally Lipschitz continuous and $J=\infty$ and $K=\infty$ hold. Then for any $u_0 \in C({\BBR}^N)$ which satisfies
	\bel {Q1}0 \leq u_0(x) \leq w_b(\abs x) \q \forall x \in {\BBR}^N	\ee
for some $0<b$, there exists a positive solution $\unl u \in C(\ovl {Q_\ity})$ of $(\ref{A1})$ in $Q_\ity$ satisfying $\unl u(.,0)=u_0$ in $\BBR^N$ and 
		\bel {Q2} \unl u(x,t) \leq \min\{w_b(\abs x),\gf_\ity(t)\} \forevery (x,t) \in Q_\ity. \ee
\es
\Proof For any $R>0$, let $u_R$ be the solution of
	$$ \left\{ \BA{lll} 
	\prt_t u_R - \Gd_p u_R +f(u_R) &= 0 &\text{in } Q_\ity \\
	\phantom{\prt_t u_R - \Gd u_R -}
	u_R(x,0)&=u_0(x)\gc_{B_R}(x) \qq &\text{in } {\BBR}^N.
	\EA \right. $$
The functions $\gf_\ity$ and $w_b$ are solutions of $(\ref{A1})$ in $Q_\ity$, which dominate $u_R$ at $t=0$, therefore, by the maximum principle, 
	\bel{Ex3.1} \min\{\gf_\ity(t), w_b(\abs x)\} \geq u_R(x,t) \forevery (x,t) \in Q_\ity. \ee 
The mapping $R \mt u_R$ is increasing,  jointly with $(\ref{Ex3.1})$ it implies that there exists a solution $\unl u:={\dsps \lim_{R \to \ity}u_R}$ of  $(\ref{A1})$ in $Q_\ity$ which satisfies $\unl u(x,0)=u_0(x)$ in ${\BBR}^N$. Letting $R \to \ity$ in $(\ref{Ex3.1})$ yields to $(\ref{Q2})$. \qeda  \medskip

\noindent {\bf Proof of \rth{non-unique}.} Combining \rprop{Ex2} and \rprop{Ex3} we see that there exist two solutions $\unl u$ and $\ovl u$ with the same initial data $u_0$, which are ordered and different since ${\dsps \lim_{\abs x \to \ity} \ovl u(x,t) = \ity}$ and ${\dsps \lim_{\abs x \to \ity}\unl u(x,t) \leq \gf_\ity(t) < \ity}$ for all $t >0$. \qeda
\mysection{Estimate and stability}
In this section we assume that $\Gw$ is a domain in $\BBR^N$, possibly unbounded, $0<T\leq \ity$ and set $Q_T^\Gw:=\Gw \ti (0,T)$ and $Q_T:=\BBR^N \ti (0,T)$. We denote by $\GTM(\Gw)$ the set of Radon measures in $\Gw$ and by $\GTM_+(\Gw)$ its positive cone. 
\bdef{weak-sol}
A nonnegative function $u$ is called a {\it weak solution} of $(\ref{A1})$ in $Q_T^\Gw$ if $u$, $\abs{\nabla u}^p$, $f(u) \in L_{loc}^1(Q_T^\Gw)$ and
	\bel{E1} \myint{0}{T}\myint{\Gw}{}\left(-G(u)\prt_t \vgf + \abs{\nabla u}^{p-2}{\nabla u}.\nabla(g(u)\vgf) + f(u)g(u)\vgf\right)dxdt = 0  \ee
for any $\vgf \in C_c^\ity(Q_T^\Gw)$ and any function $g \in C(\BBR) \cap W^{1,\ity}(\BBR)$ where $G'(r)=g(r)$. 
\es


The next results are obtained by adapting the proofs in \cite{BvCV}.
\subsection{Regularity Properties}
The following integral estimates are essentially \cite[Prop 2.1]{BvCV} with $u^q$ replaced by $f(u)$.
\bprop{IE} Assume $p>1$. Let $\gd <0$, $\gd \ne -1$ and $0<t<\theta<T$. Let $u$ be a nonnegative weak solution of $(\ref{A1})$ in $Q_T^\Gw$. For any nonnegative function $\zeta \in C_c^\ity(\Omega)$ and $\tau>p$,
	\bel{IEa} \BA{ll}
	\myfrac{1}{\gd+1}\myint{\Gw}{}(1+u(x,t))^{1+\gd}\zeta^\tau(x)dx + \myfrac{\abs \gd}{2}\myint{t}{\theta}\myint{\Gw}{}(1+u)^{\gd-1}\zeta^\tau\abs{\nabla u}^pdx\,dt \\ [4mm] \phantom{-}
	\leq \myfrac{1}{\gd+1}\myint{\Gw}{}(1+u(x,\theta))^{1+\gd}\zeta^\tau (x)dx + \myint{t}{\theta}\myint{\Gw}{}(1+u)^\gd f(u)\zeta^\tau dx\,dt \\ [4mm] \phantom{---------}
	+c_{18}\myint{t}{\theta}\myint{\Gw}{}(1+u)^{\gd+p-1}\zeta^{\tau-p}\abs{\nabla \zeta}^pdx\,dt.
	\EA \ee 
and
	\bel{IEb} \BA{ll}
	\myint{\Omega}{}(1+u(x,t))\zeta^\tau (x)dx \leq \myint{\Omega}{}(1+u(x,\theta))\zeta^\tau (x)dx + \myint{t}{\theta}\myint{\Omega}{}f(u)\zeta^\tau dx\,dt \\ [4mm] \phantom{-}
	+\tau\myint{t}{\theta}\myint{\Omega}{}(1+u)^{\gd-1}\zeta^\tau \abs{\nabla u}^p dx\,dt + \tau\myint{t}{\theta}\myint{\Omega}{}(1+u)^{(1-\gd)(p-1)}\zeta^{\tau-p} \abs{\nabla \zeta}^p dx\,dt.
	\EA \ee
Conversely,
	\bel{IEb4} \BA{ll}
	\myfrac{1}{4}\myint{\Gw}{}u(x,\theta)\zeta^\tau(x)dx + \myfrac{1}{2}\myint{t}{\theta}\myint{\Gw}{}f(u)\zeta^\tau dx\,dt \\ [4mm] \phantom{--}
	\leq c_{19} + \myint{\Gw}{}u(x,t)\zeta^\tau(x)dx + \tau\myint{t}{\theta}\myint{\Gw}{}\zeta^{\tau-1}\abs{\nabla u}^{p-1}\abs{\nabla \zeta}dx\,dt 
	\EA \ee
	and 
	\bel{IEc} \BA{ll}
	\myfrac{1}{4}\myint{\Omega}{}(1+u(x,\theta))\zeta^\tau (x)dx + \myfrac{1}{2}\myint{t}{\theta}\myint{\Omega}{}f(u)\zeta^\tau dx\,dt \leq \myint{\Gw}{}(1+u(x,t))\zeta^\tau (x)dx \\ [4mm] 
	+ \,\, \tau\myint{t}{\theta}\myint{\Gw}{}(1+u)^{\gd-1}\zeta^\tau \abs{\nabla u}^pdx\,dt + \tau\myint{t}{\theta}\myint{\Gw}{}(1+u)^{(1-\gd)(p-1)}\zeta^{\tau-p}\abs{\nabla \zeta}^p dx\,dt + c_{20}
	\EA \ee 
where $c_i=c_i(p,f,\gd,\tau)$ $(i=18,19,20)$.
\es

The next result is the keystone for the existence of an initial trace in the class of Radon measures. It is  essentially \cite[Prop 2.2]{BvCV} with $u^q$ replaced by $f(u)$, but we shall sketch its proof for the sake of completeness.
\bprop{RP} Let $u$ be a nonnegative solution of $(\ref{A1})$ in $Q_T^\Gw$. Let $0<\theta<T$. Assume that two of the three following conditions holds, for any open set $U \sbs \sbs \Gw$:
	\bel{RP1} {\dsps \sup_{t \in (0,\theta]}\myint{U}{}u(x,t)dx<\ity}, \ee
	\bel{RP2} \myint{0}{\theta}\myint{U}{}f(u)dx\,dt<\ity, \ee	
	\bel{RP3} \myint{0}{\theta}\myint{U}{}\abs{\nabla u}^{p-1}dx\,dt<\ity. \ee
Then the third one holds for any $U \sbs \sbs \Gw$. Moreover,
	\bel{RP4} \myint{0}{\theta}\myint{U}{}u^\sigma dx\,dt < \ity \forevery \sigma \in (0,q_c) \ee
and
	\bel{RP5} \myint{0}{\theta}\myint{U}{}\abs{\nabla u}^r dx\,dt < \ity \forevery r \in (0,\myfrac{N}{N+1}q_c) \ee
where $q_c=p-1+p/N$. Finally, there exists a Radon measure $\mu \in \GTM_+(\Gw)$ such that for any $\zeta \in C_c(\Gw)$,
	\bel{IT1} {\dsps \lim_{t \to 0}\myint{\Gw}{}u(x,t)\zeta(x)dx}=\myint{\Gw}{}\zeta(x)d\mu \ee
and $u$ satisfies
	\bel{IT2} \BA{ll} \myint{0}{\theta}\myint{\Gw}{}(-u\prt_t \vgf + \abs{\nabla u}^{p-2}\nabla u. \nabla \vgf + f(u)\vgf)dx\,dt \\
	\phantom{----------------}
	=\myint{\Gw}{}\vgf(x,0)d\gm - \myint{\Gw}{}u(x,\theta)\vgf(x,\theta)dx
	\EA \ee
for any $0<\theta<T$ and $\vgf \in C_c^\ity(\Gw\ti[0,T))$.
\es 
\Proof \noindent{\it Step 1: Assume $(\ref{RP1})$ and $(\ref{RP3})$ holds.} Let $\zeta$ and $\tau$ as in \rprop{IE}, there holds
\bel{IEb1}\BA {ll}
\myint{\Gw}{}(1+u(x,t))\gz^\gt dx=\myint{\Gw}{}(1+u(x,\gth))\gz^\gt dx+\myint{t}{\gth}\myint{\Gw}{}
f(u)\gz^\gt dxdt\\[4mm]\phantom{----------------}
+\gt\myint{t}{\gth}\myint{\Gw}{}\gz^{\gt-1}|\nabla u|^{p-2}\nabla u.\nabla\gz dxdt.
\EA\ee
It follows  that $f(u) \in L^1((0,\theta),L_{loc}^1(\Gw))$. \medskip

\noindent{\it Step 2: Assume that $(\ref{RP2})$ and $(\ref{RP3})$ hold.} Then $(\ref{RP1})$ follows from $(\ref{IEb1})$. \medskip

\noindent{\it Step 3: Assume that $(\ref{RP1})$ and $(\ref{RP2})$ hold.} Let $\gd \in (\max(1-p,-1),0)$ be fixed. From $(\ref{IEa})$, we get for any $0<t<\theta$,
	\bel{RP3a} \BA{l}\myfrac{\abs \gd}{2}\myint{t}{\theta}\myint{\Gw}{}(1+u)^{\gd-1}\abs{\nabla u}^p\zeta^\tau dx\,dt \leq \myfrac{1}{\gd+1}\myint{\Gw}{}(1+u(x,\theta))^{\gd+1}\zeta^\tau dx 
	\\ [4mm] \phantom{----} 
	+ \myint{t}{\theta}\myint{\Gw}{}(1+u)^\gd f(u)\zeta^\tau dx\,dt 
	+ c_{18}\myint{t}{\theta}\myint{\Gw}{}(1+u)^{\gd+p-1}\zeta^{\tau-p}\abs{\nabla \zeta}^pdx\,dt. \EA \ee
If $p \leq 2$, then $(1+u)^{\gd +p-1} \leq 1+u$. Consequently, by $(\ref{RP1})$,
	\be \myint{0}{\theta}\myint{\Gw}{}(1+u)^{\gd+p-1}\zeta^{\tau-p}\abs{\nabla \zeta}^pdx\,dt < \myint{0}{\theta}\myint{\Gw}{}(1+u)\zeta^{\tau-p}\abs{\nabla \zeta}^pdx\,dt<\ity, \ee
which, along with $(\ref{RP2})$ and $(\ref{RP3a})$, implies that
	\bel{RP3b} \myint{t}{\theta}\myint{\Gw}{}(1+u)^{\gd-1}\abs{\nabla u}^p\zeta^\tau dx\,dt < c_{21}. \ee
If $p>2$, we choose $\gd \in (1-p,2-p)$, $\gd \ne -1$, $\zeta$ and $\tau$ as in \rprop{IE}, then $(\ref{IEa})$ remains valid. From the inequality $(1+u)^{1+\gd} < 1+u$ and $(\ref{RP1})$, we find that 
	$$ \myfrac{1}{\abs{\gd+1}}\myint{\Gw}{}(1+u(x,t))^{1+\gd}\zeta^\tau(x)dx < \myfrac{1}{\abs{\gd+1}}\myint{\Gw}{}(1+u(x,t))\zeta^\tau(x)dx< c_{22}. $$
Hence, by $(\ref{IEa})$, 
	\bel{RP3a*} \BA{l}\myfrac{\abs \gd}{2}\myint{t}{\theta}\myint{\Gw}{}(1+u)^{\gd-1}\zeta^\tau \abs{\nabla u}^p dx\,dt \leq \myfrac{1}{\gd+1}\myint{\Gw}{}(1+u(x,\theta))^{\gd+1}\zeta^\tau dx 
	\\ [4mm] \phantom{--} 
	+ \myint{t}{\theta}\myint{\Gw}{}(1+u)^\gd f(u)\zeta^\tau dx\,dt 
	+ c_{18}\myint{t}{\theta}\myint{\Gw}{}(1+u)^{\gd+p-1}\zeta^{\tau-p}\abs{\nabla \zeta}^pdx\,dt + c_{22}. \EA \ee 
Since $\gd<2-p$, $\gd+p-1<1$, hence $(1+u)^{\gd+p-1}\leq 1+u$. Therefore, $(\ref{RP3b})$ follows from $(\ref{RP1})$, $(\ref{RP2})$ and $(\ref{RP3a*})$.

By applying the Gagliardo-Nirenberg inequality as in \cite[Prop 2.2 (iii)]{BvCV}, we deduce that
	$$ \myint{0}{\theta}\myint{U}{}(1+u(x,t))^\sigma dx < c_{23}$$
for any $\sigma \in (0,q_c)$ with $q_c=p-1+p/N$, which leads to $(\ref{RP4})$. Next for $0<r<p$, and any $\gd <0$, we find
	\bel{RP3g} \BA{l} \myint{0}{\theta}\myint{U}{}\abs{\nabla u}^r dx \leq \bigg( \myint{0}{\theta}\myint{U}{}(1+u)^{\gd-1}\abs{\nabla u}^p dx\,dt\bigg)^\frac{r}{p} \\ [3mm] \phantom{--------}
	\ti \bigg( \myint{0}{\theta}\myint{U}{}(1+u)^\frac{(1-\gd)r}{p-r} dx\,dt\bigg)^\frac{p-r}{p}. \EA \ee
Thus, if $r \in (0,Nq_c/(N+1))$, this proves $(\ref{RP5})$; furthermore, since $p-1<Nq_c/(N+1)$, we obtain $(\ref{RP3})$. \medskip

\noindent{\it Step 4: End of the proof.} Now we use $(\ref{E1})$ with $g=1$, for any $\gz \in C_c^\ity(\Gw)$ and any $0<t<\theta<T$,
	\bel{RP6}  \myint{\Gw}{}u(x,t)\gz(x)dx = \myint{\Gw}{}u(x,\theta)\gz(x)dx + \myint{t}{\theta}\myint{\Gw}{}\left(\abs{\nabla u}^{p-2}\nabla u.\nabla \gz + f(u)\gz\right)dx\,dt. \ee
Because the right-hand side of $(\ref{RP6})$ has a finite limit when $t \to 0$, the same holds with $t \mapsto \myint{\Gw}{}u(x,t)\gz(x)dx$. The mapping $\gz \mapsto \lim_{t \to 0}\myint{\Gw}{}u(x,t)\gz(x)dx$ is a positive linear functional $\ell_\Gw$ on the space $C_c^\ity(\Gw)$. By a partition of unity it can be extended in a unique way as a Radon measure $\mu \in \GTM_+(\Gw)$ and $(\ref{IT1})$ holds.  \smallskip

Finally, let $0<t<\theta$ be fixed,  $g=1$ and $\vgf \in C_c^\ity(Q_T^\Gw)$, thus
	 \bel{RP7} \BA{ll} \myint{t}{\theta}\myint{\Gw}{}(-u\prt_t \vgf + \abs{\nabla u}^{p-2}\nabla u. \nabla \vgf + f(u)\vgf)dx\,d\tau \\
	\phantom{-------}
	=\myint{\Gw}{}u(x,t)\vgf(x,0)dx - \myint{\Gw}{}u(x,\theta)\vgf(x,\theta)dx.
	\EA \ee
But
	$$ \abs{\myint{\Gw}{}u(x,t)(\vgf(x,t)-\vgf(x,0))dx}\leq c_{24}t\myint{\Gw}{}u(x,t)dx. $$
By $(\ref{IT1})$, letting $t \to 0$ yields
	$$ \myint{\Gw}{}u(x,t)\vgf(x,t)dx \to \myint{\Gw}{}\vgf(x,0)d\gm. $$
Thus, letting $t \to 0$ in $(\ref{RP7})$ implies $(\ref{IT2})$. \qeda \medskip 

Next we consider the the following problems
	\bel{Pu} \left\{ \BA{lll} \prt_t u - \Gd_p u + f(u)&= 0 \qq &\text{in } Q_T^\Gw, 
	\\	\phantom{-------,}
	u&=0 &\text{in } \prt\Gw\ti (0,T)
	\\	\phantom{-----,}
	u(.,0)&=\gm &\text{in } \Gw.
	\EA \right. \ee
where $\gm \in \GTM_+(\Gw)$. The solutions are considered in the entropy sense (see \cite{SLT} and \cite{Li}).\medskip

We recall that for $q\geq 1$ and $\Gth\subset\BBR^d$ open, the Marcinkiewicz space (or weak Lebesgue space) $M^q(\Gth)$ is the set of all locally integrable functions $u:\Gth\mapsto\BBR$ such there exists $C\geq 0$ with the property that for any measurable set $E\subset \Gth$,
\bel{Mar}\myint{E}{}|u|dy\leq C|E|^{1-\frac{1}{q}}.
\ee
The norm of $u$ in $M^q(\Gth)$ is the smallest constant such that $(\ref{Mar})$ holds for any measurable set $E$ (see \cite{SLT}, \cite{Li} for more details). Here $dy$ denotes the Lebesgue measure in $\BBR^d$, although any positive Borel measure can be used.  

We recall the following result of Segura de Leon and Toledo \cite[Th 2]{SLT} and Li \cite[Th 1.1]{Li} dealing with entropy solutions with initial data in $L^1$. However such solutions coincide with the semi-group solutions because of uniqueness.
\bprop{reg-Mar} Assume $p>\frac{2N}{N+1}$, $\Gw\subset \BBR^N$ is any open subset  and, $h\in L^1(Q^\Gw_T )$and  $\gm \in L_+^1(\Gw)$. Let $v\in C([0,T;L^1(\Gw))$ be the entropy solution to problem
	\bel{Pv} \left\{ \BA{lll} \prt_t v - \Gd_p v &\!\!= h \qq &\text{in } Q^\Gw_T \\
	\phantom{----,}
	v&\!\!=0&\text{in } \prt \Gw\ti (0,\infty)\\
	\phantom{--,}
	v(.,0)&\!\!=\gm &\text{in } \Gw.
	\EA \right. \ee
Then $v \in M^{p-1+\frac{p}{N}}(Q^\Gw_T)$, $\nabla v\in M^{p-\frac{N}{N+1}}(Q^\Gw_T)$ and there holds
\bel{Mar2}\norm v_{M^{p-1+\frac{p}{N}}(Q^\Gw_T)}+\norm {\nabla v}_{M^{p-\frac{N}{N+1}}(Q^\Gw_T)}\leq c_{25},
\ee
for some $c_{25}>0$ depending on $p$, $N$, $\norm \gm_{L^1(\Gw)}$ and $\norm h_{L^1(Q^\Gw_T)}$.
\es


\subsection{Stability}
Let $\{\gm_n\}\sbs L_+^1(\BBR^N)$ be a sequence converging to $\gm$ in weak sense of measures, then $\norm{\mu_n}_{L^1(\BBR^N)} \leq c^*$, where $c^*$ depends only on $N,p$ and $\norm{\gm}_{\GTM(\BBR^N)}$. Denote by $u_{\gm_n}$ (resp. $v_{\gm_n}$) the solution to problem $(\ref{Pu})$ (resp. $(\ref{Pv})$ with $h \equiv 0$) with the initial data $\mu_n$. Then the following estimate holds
	\bel{stab1} 0 \leq u_{\gm_n} \leq v_{\gm_n}. \ee
By \cite[Theorem 3]{HeVa}, 
	$$ \norm{v_{\gm_n}(.,t)}_{L^\ity(\BBR^N)} \leq c_{26}t^\frac{-N}{N(p-2)+p}\norm{\gm_n}_{L^1(\BBR^N)}^\frac{p}{N(p-2)+p} \forevery t>0, $$
where $c_{26}=c_{26}(N,p)>0$. Thus 
 	\bel{stab2} \BA{ll} 
	\norm{u_{\gm_n}(.,t)}_{L^\ity(\BBR^N)} \leq c_{26}t^\frac{-N}{N(p-2)+p}\norm{\gm_n}_{L^1(\BBR^N)}^\frac{p}{N(p-2)+p} \\[3mm] 
	\phantom{-------,,}
	\leq c_{27} t^\frac{-N}{N(p-2)+p}
	\EA \ee
for every $t>0$, where $c_{27}=c_{27}(N,p,c^*)>0$. 

It follows from $(\ref{Mar2})$ and $(\ref{stab1})$ that
 	\bel{stab3}  \norm{u_{\gm_n}}_{M^{p-1+p/N}(Q_T)} \leq c_{25}\norm{\gm_n}_{L^1(\BBR^N)}^\frac{p+N}{1+p(N-1)} \leq c_{28}(N,p,c^*). \ee
By $(\ref{stab2})$ and the regularity theory of degenerate parabolic equations \cite{Di}, we derive that the sequence $\{u_{\gm_n}\}$ is equicontinuous in any compact subset of $Q_T$. As a consequence, there exist a subsequence, still denoted by $\{u_{\gm_n}\}$ and a function $u$ such that $\{u_{\gm_n}\}$ converges to $u$ locally uniformly in $Q_T$. 

\blemma{stab-lem1} The sequence $f(u_{\gm_n})$ converges strongly to $f(u)$ in $L^1(Q_T)$. Furthermore, $\{u_n\}$ converges strongly to $u$ in $L_{loc}^q(\overline{Q_T})$ for every $1\leq q< q_c$.
\es
\Proof Since $u_{\gm} \to u$ a.e in $Q_T$, by Vitali's theorem, it is sufficient to show that the sequence $\{f(u_{\gm_n})\}$ is uniformly integrable.
Let $E$ be a Borel subset of $Q_T$ and let $R>0$. Then, since $f$ is increasing,
	$$ \BA{ll} \myint{}{}\myint{E}{}f(u_{\gm_n})dx\,dt =\myint{}{}\myint{E \cap \{u_{\gm_n\leq R}\}}{}f(u_n)dx\,dt +  \myint{}{}\myint{E \cap \{u_{\gm_n}> R\}}{}f(u_{\gm_n})dx\,dt \\[4mm]
	\phantom{--------} 
	\leq f(R)\myint{}{}\myint{E}{}dx\,dt+\myint{}{}\myint{E \cap \{u_{\gm_n}> R\}}{}f(u_{\gm_n})dx\,dt.
	\EA $$
For $\gl \geq 0$, we set $B_n(\gl)=\{(x,t )\in Q_T): u_{\gm_n}>\gl\}$ and $a_n(\gl)=\myint{}{}\myint{B_n(\gl)}{}dx\,dt$. Then
 \bel{ex-est4} \myint{}{}\myint{E \cap \{u_{\gm_n}> R\}}{}f(u_{\gm_n})dx\,dt \leq \myint{}{}\myint{\{u_{\gm_n}\geq R\}}{}f(u_{\gm_n})dx\,dt = - \myint{R}{\ity}f(\gl)da_n(\gl) \ee
and
	$$ -\myint{R}{\ity}f(\gl)da_n(\gl) \leq f(R)a_n(R) + \myint{R}{\ity}a_n(\gl)df(\gl). $$
It follows from $(\ref{stab3})$ that
	$$ a_n(\gl) \leq c_{25}\norm{\gm_n}_{\GTM_+(\BBR^N)}^\frac{p+N}{1+p(N-1)}\gl^{-(p-1+\frac{p}{N})}\leq c_{29}\gl^{-(p-1+\frac{p}{N})}. $$
Plugging these estimates into $(\ref{ex-est4})$ yields
 \bel{ex-est5} \BA{ll} 
  \myint{}{}\myint{E \cap \{u_{\gm_n}> R\}}{}f(u_{\gm_n})dx\,dt \leq f(R)a_n(R)+c_{29}\myint{R}{\ity}\gl^{-(p-1+\frac{p}{N})}df(\gl) \\
  \phantom{\myint{}{}\myint{E \cap \{u_{\gm_n}> R\}}{}f(u_{\gm_n})dx\,dt}
  \leq f(R)a_n(R)-c_{29}f(R)R^{-(p-1-\frac{p}{N})}\\
  \phantom{\myint{}{}\myint{E \cap \{u_{\gm_n}> R\}}{}f(u_{\gm_n})dx\,dt}
  +c_{29}\left(p-1+\myfrac{p}{N}\right)\myint{R}{\ity}f(\gl)\gl^{-(p+\frac{p}{N})}d\gl \\
  \phantom{\myint{}{}\myint{E \cap \{u_{\gm_n}> R\}}{}f(u_{\gm_n})dx\,dt}
  \leq c_{29}\left(p-1+\myfrac{p}{N}\right)\myint{R}{\ity}f(\gl)\gl^{-(p+\frac{p}{N})}d\gl.
  \EA \ee
Since 
	$$ \myint{1}{\ity}\gl^{-(p+\frac{p}{N})}f(\gl)d\gl < \ity, $$
for given $\ge>0$,we can choose $R>0$ large enough such that
	$$ c_{29}\left(p-1+\myfrac{p}{N}\right)\myint{R}{\ity}f(\gl)\gl^{-(p+\frac{p}{N})}d\gl<\ge/2. $$
Set $\gd=(1+f(R))^{-1}\ge/2$, then
	$$ \abs{E} <\gd \Lra 0 \leq \myint{}{}\myint{E}{}f(u_{\gm_n})dx\,dt < \ge, $$
which proves the uniform integrability of the sequence $\{f(u_{\gm_n})\}$. The last assertion follows from the fact that $u_{\gm_n}$ is bounded in $M^{q_c}(Q_T)$ (remember that $q_c=p-1+p/N$) and $M^{q_c}(Q_T)\subset L^q_{loc}(\overline{Q_T})$ with continuous imbedding, for any $q<q_c$. The conclusion follows again by Vitali's theorem.
\qeda

\blemma{stab-lem2} Assume $p>\frac{2N}{N+1}$, then for any $U \sbs \sbs \BBR^N$, the sequence $\{\nabla u_{\gm_n}\}$ converges strongly to $\nabla u$ in $(L^{s}(Q_T))^N$ for every $1\leq s<s_c:=p-\frac{N}{N+1}$. \es
\Proof We set $h_n=-f(u_{\gm_n})$ and write the equation under the form
	\bel{Pv2} \left\{ \BA{lll} \prt_t u_{\gm_n} - \Gd_p u_{\gm_n} &\!\!\!\!= h_n \qq &\text{in } Q_T \\
	\phantom{--.,,}
	u_{\gm_n}(.,0)&\!\!\!\!=\gm_n &\text{in } \BBR^N.
	\EA \right. \ee
We already know from the $L^1$-contraction principle and \rprop {reg-Mar} that 
$$\norm{u_{\gm_n}(.,t)}_{L^1(\BBR^N)}\leq \norm{\gm_n}_{L^1(\BBR^N)}\qquad\forall t\in (0,T]
$$
and $u_{\gm_n}\to u$ in $L^q_{loc}(\overline{Q_T})$ for every $q\in [1,q_c)$ and $|\nabla u_{\gm_n}|$ is bounded in $L^s_{loc}(\overline{Q_T})$ for every $1\leq s<s_c$. Thus $|\nabla u_{\gm_n}|^{p-1}$ remains bounded in bounded in $L^\gs_{loc}(\overline{Q_T})$ for every $1\leq \gs<\gs_c:=1+\frac{1}{(N+1)(p-1)}$. Furthermore,  
\bel{Cau-D} \{\nabla u_{\gm_n}\} \text{ is a Cauchy sequence in measure}. \ee
and the proof is similar to  the one of \cite[Th 5.1-step2]{BvCV}. Up to the extraction of a subsequence, $\{\nabla u_{\gm_n}\}$ converges a.e. to some $D=(D_1,...,D_N)$ in $Q_T$. Consequently, $\{\abs{\nabla u_{\gm_n}}^{p-2}\nabla u_{\gm_n}\}$ converges a.e. to $\abs{D}^{p-2}D$ in $Q_T$ and, by Vitali's theorem,
	\bel{converg3} \BA {cll}
	\nabla u_{\gm_n}\to D \qq &\text{strongly in } (L^s_{loc}(\overline{Q_T}))^N,\q \forall s \in [1,s_c),
	\\[2mm]
	\{\abs{\nabla u_{\gm_n}}^{p-2}\nabla u_{\gm_n}\} \to \abs{D}^{p-2}D \qq &\text{strongly in } (L^\gs_{loc}(\overline{Q_T}))^N,\q \forall \gs \in [1,\gs_c). 	\EA	\ee
	 which implies $\nabla u = D$ and the conclusion of the lemma follows. \qeda\medskip

\noindent{\bf Proof of \rth{stab}.} 
{\it Step 1.} For any $\zeta \in C_c^\ity(\BBR^N)$ and $t>0$, we have
	$$ \myint{\BBR^N}{}u_{\gm_n}(x,t)\zeta(x)dx+\myint{0}{t}\myint{\BBR^N}{}(\abs{\nabla u_{\gm_n}}^{p-2}\nabla u_{\gm_n}\nabla \zeta+f(u_{\gm_n})\zeta)dx\,dt  =\myint{\BBR^N}{}\gm_n(x)\zeta(x)dx$$
By \rlemma{stab-lem1} and \rlemma{stab-lem2}, up to the extraction of a subsequence, we can pass to the limit in each term and get
	$$ \myint{\BBR^N}{}u(x,t)\zeta(x)dx+\myint{0}{t}\myint{\BBR^N}{}(\abs{\nabla u}^{p-2}\nabla u \nabla \zeta 	+ f(u)\zeta)dx\,dt=\myint{\BBR^N}{}\zeta d\gm. $$
Letting $t \to 0$ yields
	\bel{stab4} \lim_{t \to 0}\myint{\BBR^N}{}u(x,t)\zeta(x)dx=\myint{\BBR^N}{}\zeta(x). \ee
For any $\vgf \in C_c^\ity(\BBR^N \ti [0,\ity))$ and $\theta>0$, we have
	\bel{stab5} \BA{ll} \myint{0}{\theta}\myint{\BBR^N}{}(-u_{\gm_n}\prt_t\vgf + \abs{\nabla u_{\gm_n}}^{p-2}\nabla u_{\gm_n}.\nabla \vgf + 		
	f(u_{\gm_n})\vgf)dx\,dt \\
	\phantom{---}
	=\myint{\BBR^N}{}\vgf(0,x)\gm_n(x)dx-\myint{\BBR^N}{}u_{\gm_n}(x,\theta)\vgf(x,\theta)dx. 
	\EA \ee
By the previous convergence results, we can pass to the limit in $(\ref{stab5})$ to obtain
	\bel{stab8} \BA{ll} \myint{0}{\theta}\myint{\BBR^N}{}(-u\prt_t\vgf + \abs{\nabla u}^{p-2}{\nabla u\nabla \vgf + f(u)\vgf})dx\,dt \\ 
	\phantom{-----------}
	= \myint{\BBR^N}{} \vgf(.,0)d\gm - \myint{\BBR^N}{}u(.,\theta)\vgf(.,\theta)dx. \EA \ee
\noindent{\it Step 2: $u$ is a weak solution.} By $(\ref{stab2})$ 
$$\sup\{\norm{u_{\gm_n}(.,t)}_{L^\infty(\BBR^N)} ,\norm{u(.,t)}_{L^\infty(\BBR^N)}\}\leq c_{27}t^{-\frac{N}{N(p-2)+p}}\qquad\forall t\in (0,T].
$$
 Let $\zeta \in C_c^\ity(\BBR^N)$. Since $\{u_{\gm_n}(.,\theta)\}$ converges locally uniformly to $u(.,\theta)$ in $\BBR^N$, for any $\gth>0$, there holds 
	\bel{stab10} \BA{ll} \myfrac{1}{2}\myint{\BBR^N}{}(u_{\gm_n}-u_{\gm_m})^2(.,T)\zeta dx\,dt + 	
	\myint{\theta}{T}\myint{\BBR^N}{}(f(u_{\gm_n})-f(u_{\gm_m}))(u_{\gm_n}-u_{\gm_m})\zeta dx\,dt \\[4mm]
	\phantom{---}
	+ \myint{\theta}{T}\myint{\BBR^N}{}(\abs{\nabla u_{\gm_n}}^{p-2}\nabla u_{\gm_n}-\abs{\nabla u_{\gm_m}}^{p-2}\nabla 	
	u_{\gm_m}).\nabla (u_{\gm_m}- u_{\gm_n})\zeta dx\,dt \\[4mm]
	\phantom{---}
	\leq \myfrac{1}{2}\myint{\BBR^N}{}(u_{\gm_n}-u_{\gm_m})^2(.,\theta)\zeta dx\,dt \\[4mm]
	\phantom{---}
	+ \myint{\theta}{T}\myint{\BBR^N}{}\abs{\abs{\nabla 
	u_{\gm_n}}^{p-2}\nabla u_{\gm_n}-\abs{\nabla u_{\gm_m}}^{p-2}\nabla 	
	u_{\gm_m}}\abs{u_{\gm_m}-u_{\gm_n}}\abs{\nabla \zeta} dx\,dt.
	\EA \ee
This implies directly
	\bel{stab11} \nabla u_{\gm_n} \to \nabla u \text{ in } L^p_{loc}(Q_T), \ee
by \rlemma{stab-lem2} when $p\geq 2$. When $1<p<2$, we derive by Fatou's lemma
	\bel{stab11'} \BA{ll} \myfrac{1}{2}\myint{\BBR^N}{}(u_{\gm_n}-u)^2(.,T)\zeta dx\,dt + 	
	\myint{\theta}{T}\myint{\BBR^N}{}(f(u_{\gm_n})-f(u))(u_{\gm_n}-u)\zeta dx\,dt \\[4mm]
	\phantom{---}
	+ \myint{\theta}{T}\myint{\BBR^N}{}(\abs{\nabla u_{\gm_n}}^{p-2}\nabla u_{\gm_n}-\abs{\nabla u}^{p-2}\nabla 	
	u).\nabla (u_{\gm_n}- u)\zeta dx\,dt \\[4mm]
	\phantom{---}
	\leq \myfrac{1}{2}\myint{\BBR^N}{}(u_{\gm_n}-u)^2(.,\theta)\zeta dx\,dt \\[4mm]
	\phantom{---}
	+ \myint{\theta}{T}\myint{\BBR^N}{}\abs{\abs{\nabla 
	u_{\gm_n}}^{p-2}\nabla u_{\gm_n}-\abs{\nabla u}^{p-2}\nabla 	
	u}\abs{u_{\gm_m}-u}\abs{\nabla \zeta} dx\,dt.
	\EA \ee
Using again \rlemma{stab-lem2}, it implies
\bel{stab11''} \lim_{n\to\infty}\myint{\theta}{T}\myint{\BBR^N}{}\abs{\nabla u_{\gm_n}}^{p}\gz dx\,dt=
\myint{\theta}{T}\myint{\BBR^N}{}\abs{\nabla u}^{p}\gz dx\,dt.
 \ee
Since $\nabla u_{\gm_n}\rightharpoonup \nabla u$ weakly in $L^p_{loc}(Q_T)$, it implies again that 
$(\ref{stab11})$ holds true.
At end,  let $\vgf \in C_c^\ity(Q_T)$ and consider $0<\theta<T$ and $U \sbs\sbs \BBR^N$ such that $\supp\vgf \sbs (\theta,T)\ti U$. Let $g \in C(\BBR^N) \cap W^{1,\ity}(\BBR^N)$ where $G'(r)=g(r)$. Multiplying the equation in $(\ref{Pu})$ (with initial data $\gm=\gm_n$) by $g(u_{\gm_n})\vgf$, we obtain
	\bel{stab12} \BA{ll} \myint{0}{T}\myint{\BBR^N}{}(-G(u_{\gm_n})\prt_t\vgf+\abs{\nabla u_{\gm_n}}^p g'(u_{\gm_n})\vgf )dx\,dt \\[3mm]
	\phantom{------}
	+ g(u_{\gm_n})\abs{\nabla u_{\gm_n}}^{p-2}\nabla u_{\gm_n}.\nabla \vgf 		
	+ \myint{0}{T}\myint{\BBR^N}{}g(u_{\gm_n})f(u_{\gm_n})dx\,dt = 0. 
	\EA \ee
By \rlemma{stab-lem1} and $(\ref{stab11})$, we can pass to the limit in each term. As a consequence, $u$ is a weak solution. \medskip

\noindent{\it Step 3: Stability.} Assume that $\{\gm_n\}$ is a sequence of functions in $L_+^1(\BBR^N)$ with compact support, which converges to $\gm \in \GTM_+^b(\BBR^N)$ in the dual sense of $C(\BBR^N)$, then $\norm{\gm_n}_{L^1(\BBR^N)}$ is bounded independently of $n$. By the same argument as in step 1 and step 2, we can pass to the limit in each term of $(\ref{stab12})$, hence the conclusion follows. \qeda 

\blemma{sin-grad} Assume $p>2$. Let $u\in C(Q_T)$ be a positive weak solution of $(\ref{A1})$ in $Q_T$. Assume that there exists $r>0$ such that
	\bel{sin-grad1} \myint{0}{T}\myint{B_r}{}\abs{\nabla u}^{p-1}dx\,dt=\ity. \ee
Then
	\bel{sin-grad2} \sup_{\tau \in (0,T)}\myint{B_{8r}}{}u(x,\tau)=\ity. \ee
\es
\Proof By contradiction we assume that $(\ref{sin-grad2})$ does not hold. Then there exist $A_1>0$ such that
	\bel{sin-grad3} \sup_{\tau \in (0,T)}\myint{B_{8r}}{}u(x,\tau)=A_1. \ee
	
\noindent{\it Step 1: We claim that} $$ u \in L^\ity(Q_T^{B_{2r}}).$$
Since $u$ is a positive subsolution of the equation in $(\ref{Y1})$, by \cite[Theorem 4.2, Chapter V]{Di}, there exists a constant $c_{30}=c_{30}(N,p)$ such that for every $x_0 \in \BBR^N$, $0<\theta\leq t_0 <T$ and $\gs \in (0,1)$, there holds
	\bel{sin-grad5} \sup_{K_{\gs\gr} \ti (t_0-\gs\theta,t_0)}u \leq \myfrac{c_{30}\theta^\frac{1}{2}}{\gr^\frac{p}{2}(1-\gs)^\frac{N(p+1)+p}{2}}\bigg(\sup_{0<\tau<t}\abs{K_\gr}^{-1}\myint{K_\gr}{}u(x,\tau)dx\bigg)^\frac{p}{2},
	\ee
where $K_\gr(x_0)$ is the cube centered at $x_0$ and wedge $2\gr$, i.e.,
	$$ K_\gr(x_0)=\{x \in \BBR^N: \max_{1\leq i \leq N}\abs{x^i-x_0^i}<\gr\}. $$
We choose $x_0=0$, $t_0=\theta=t$, $\gs=1/2$ and $\gr=4r$, then $(\ref{sin-grad5})$ becomes
	\bel{sin-grad6} \sup_{K_{2r} \ti (\frac{t}{2},t)}u \leq 2^\frac{N(p+1)+p}{2}c_{30}t^\frac{1}{2}(4r)^\frac{-p}{2}\bigg(\sup_{0<\tau<t}\abs{K_{4r}}^{-1}\myint{K_{4r}}{}u(x,\tau)dx\bigg)^\frac{p}{2}.
	\ee
Since $B_{2r} \sbs K_{2r}$ and $K_{4r} \sbs B_{8r}$, from $(\ref{sin-grad3})$ and $(\ref{sin-grad6})$, we obtain that
	\bel{sin-grad7} \sup_{B_{2r} \ti (0,T)}u \leq 2^\frac{N-p(2N+1)}{2}c_{30}T^\frac{1}{2}r^\frac{-p(N+1)}{2}A_1^\frac{p}{2}=:A_2, 
	\ee
which implies the claim. \medskip

\noindent{\it Step 2: Let $\zeta \in C_c^\ity(\BBR^N)$ such that $\zeta \geq 0$ in $\BBR^N$, $\zeta=1$ in $B_r$ and $\abs{\nabla \zeta}\leq 1/r$. We show that}
	\bel{sin-grad8} \BA{ll} J_1(t):=\myint{0}{t}\myint{B_{2r}}{}(u+1)^\frac{-2}{p}\zeta^p\abs{\nabla u}^pdx\,d\tau < \ity, \\
	J_2(t):=\myint{0}{t}\myint{B_{2r}}{}(u+1)^\frac{2(p-1)}{p}dx\,d\tau < \ity. 	
	\EA \ee
Multiplying $(\ref{A1})$ by $(u+1)^\frac{p-2}{p}\zeta^p$ and then integrating on $\BBR^N \ti [\ge,t]$ with $0<\ge<t$, we get
	$$ \BA{ll} \myfrac{p}{2(p-1)}\myint{B_{2r}}{}(u(x,t)+1)^\frac{2(p-1)}{p}\zeta^pdx + \myfrac{p-2}{p}\myint{\ge}{t}\myint{B_{2r}}{}(u+1)^\frac{-2}{p}\zeta^p\abs{\nabla u}^pdx\,d\tau \\[4mm]
	\phantom{--------------------}	
	+ \myint{\ge}{t}\myint{B_{2r}}{}(u+1)^\frac{p-2}{p}f(u)\zeta^p dx\,d\tau \\
	\phantom{---}
	=\myfrac{p}{2(p-1)}\myint{B_{2r}}{}(u(x,\ge)+1)^\frac{2(p-1)}{p}\zeta^pdx \\ [4mm]
	\phantom{-------}
	- p \myint{\ge}{t}\myint{B_{2r}}{}(u+1)^\frac{p-2}{p}\zeta^{p-1}\abs{\nabla u}^{p-2}\nabla u \nabla \zeta dx\,d\tau,
	\EA $$
which implies that
	\bel{sin-grad9} \BA{ll} \myfrac{p-2}{p}\myint{\ge}{t}\myint{B_{2r}}{}(u+1)^\frac{-2}{p}\zeta^p\abs{\nabla u}^pdx\,d\tau \\
	\phantom{---}
	\leq \myfrac{p}{2(p-1)}\myint{B_{2r}}{}(u(x,\ge)+1)^\frac{2(p-1)}{p}\zeta^pdx \\
	\phantom{--------}
	- p \myint{\ge}{t}\myint{B_{2r}}{}(u+1)^\frac{p-2}{p}\zeta^{p-1}\abs{\nabla u}^{p-2}\nabla u \nabla \zeta dx\,d\tau.
	\EA \ee	
By Young's inequality, 
	\bel{sin-grad10} \BA{ll} p \myint{\ge}{t}\myint{B_{2r}}{}(u+1)^\frac{p-2}{p}\zeta^{p-1}\abs{\nabla u}^{p-1}\abs{ \nabla \zeta} dx\,dt\tau \\
	\phantom{---}
	\leq \myfrac{p-2}{2p}\myint{\ge}{t}\myint{B_{2r}}{}(u+1)^\frac{-2}{p}\zeta^p\abs{\nabla u}^p dx\,d\tau \\
	\phantom{--------}
	+ p\bigg(\myfrac{2p^2}{p-2}\bigg)^{p-1}\myint{\ge}{t}\myint{B_{2r}}{}(u+1)^\frac{p^2-2}{p}\abs{\nabla \zeta}^p dx\,d\tau.
	\EA \ee
It follows from $(\ref{sin-grad9})$ and $(\ref{sin-grad10})$ that
	\bel{sin-grad11} \BA{ll} \myfrac{p-2}{2p}\myint{\ge}{t}\myint{B_{2r}}{}(u+1)^\frac{-2}{p}\zeta^p\abs{\nabla u}^pdx\,d\tau \\
	\phantom{---}
	\leq \myfrac{p}{2(p-1)}\myint{B_{2r}}{}(u(x,\ge)+1)^\frac{2(p-1)}{p}\zeta^pdx \\
	\phantom{--------}
	+ p\bigg(\myfrac{2p^2}{p-2}\bigg)^{p-1}\myint{0}{t}\myint{B_{2r}}{}(u+1)^\frac{p^2-2}{p}\abs{\nabla \zeta}^p dx\,d\tau.
	\EA \ee
By $(\ref{sin-grad7})$, 
	$$ \sup_{\ge \in (0,T)}\myint{B_{2r}}{}(u(x,\ge)+1)^\frac{2(p-1)}{p}\zeta^pdx \leq c_{31}(N,p,r,\zeta,A_2) $$
and 
	$$ \myint{0}{t}\myint{B_{2r}}{}(u+1)^\frac{p^2-2}{p}\abs{\nabla \zeta}^p dx\,d\tau \leq r^{-p}\myint{0}{t}\myint{B_{2r}}{}(u+1)^\frac{p^2-2}{p} dx\,dt \leq c_{32}(N,p,r,T,A_2).
	$$
Combining the previous two estimates with $(\ref{sin-grad11})$ yields
	\bel{sin-grad12} J_1(t) \leq c_{33}(N,p,r,T,\zeta), \forevery t \in (0,T). \ee
By $(\ref{sin-grad7})$, we also find that
	\bel{sin-grad13} J_2(t) \leq c_{34}(N,p,r,T,A_2). \ee

\noindent{\it Step 3: End of proof.}
By H\"older's inequality, we get
	$$ \myint{0}{t}\myint{B_{2r}}{}\abs{\nabla u}^{p-1}\zeta^{p-1}dx\,d\tau \leq (J_1(t))^\frac{p-1}{p}(J_2(t))^\frac{1}{p}. $$
By step 2, we deduce that
	\bel{sin-grad14} \myint{0}{T}\myint{B_{2r}}{}\abs{\nabla u}^{p-1}\zeta^{p-1}dx\,dt < c_{35}(N,p,r,T,\zeta), \ee
which contradicts $(\ref{sin-grad1})$. \qeda

\mysection{Initial trace}

\subsection{The dichotomy theorem}
The dichotomy result \rth{dichotomy} is a consequence of \rprop{RP} and \rlemma{sin-grad}.  \medskip

\noindent{\bf Proof of \rth{dichotomy}} By translation we may suppose that $y=0$.\smallskip

\noindent{\it Case 1: there exists an open neighborhood $U$ of $0$ such that $(\ref{RP2})$ and $(\ref{RP3})$ hold true. } Then the statement (ii) follows from \rprop{RP}. \medskip

\noindent{\it Case 2: for any open neighborhood $U$ of $0$, $(\ref{RP2})$ or $(\ref{RP3})$ does not holds. } We first suppose that $(\ref{RP3})$ does not hold. We can choose $r>0$ such that $B_{8r} \sbs U$ and $(\ref{sin-grad1})$ holds. Then the statement (i) follows from \rlemma{sin-grad}. Suppose next that $(\ref{RP3})$ holds but $(\ref{RP2})$ does not hold, then \rprop{RP} implies that $(\ref{RP1})$ does not hold and the statement (i) follows. \qeda  

\bprop{singular-point} Assume $p>2$ and $f$ is nondecreasing and satisfies $(\ref{CFS})$. Let $u$ is a positive weak solution of $(\ref{A1})$ in $Q_\ity$ with initial trace $(\CS,\gm)$. Then for every $y \in \CS$,
\bel{singular-point1} \underline U_y(x,t):= \underline U(x-y,t) \leq u(x,t) \ee
in $Q_\ity$. 
\es
\Proof By translation we may suppose that $y=0$. Since $0 \in \CS(u)$, for any $\gh>0$ small enough 
$$ \lim_{t \to 0}\myint{B_\gh}{}u(x,t)dx=\ity. $$
For $\ge>0$, denote $M_{\ge,\gh}=\myint{B_\gh}{}u(x,\ge)dx$. For any $m>m_\gh={\displaystyle \inf_{\gs>0}M_{\gs,\gh}}$ there exists $\ge=\ge(m,\gh)$ such that $m=M_{\ge,\gh}$ and ${\displaystyle \lim_{\gh \to 0}\ge(m,\gh)=0}$. Let $\tl u_{\gh}$ be the solution to the problem
$$ \left\{ \BA{lll}
\prt_t \tl u_\gh - \Gd_p \tl u_\gh + f(\tl u_\gh) &= 0 \qq &\text{in } Q_\infty \\ \phantom{\prt_t v_{\gh} - \Gd v_{}-,,}
\tl u_{\gh}(x,0)&=u(x,\ge)\gc_{_{B_\gh}} \q &\text{in } \BBR^N
\EA \right. $$
where $\gc_{_{B_\gh}}$ is the characteristic function of $B_\gh$. By the maximum principle $\tl u_{\gh}\leq u$ in $\BBR^N \ti (\ge,\ity)$. By \rth{stab} $v_{\gh}$ converges to $u_k$ when $\gh$ goes to zero. Letting $m$ go to infinity yields $(\ref{singular-point1})$. \qeda \medskip

\noindent{\bf Proof of \rth{Strong}} The conclusion follows directly from \rprop{singular-point}. \qeda

\subsection{The Keller-Osserman condition does not hold}
\blemma {Jfinite}Assume $p>2$, $(\ref{CFS})$ and $J<\infty$ are satisfied and ${\dsps \lim_{k \to \ity}u_k=\gf_\ity}$. If $u$ is a positive solution of $(\ref{A1})$ in $Q_{\infty}$ which satisfies
  	\bel{J-1}
  	\limsup_{t\to 0}\myint{G}{}u(x,t)dx=\infty,
  	\ee
for some bounded open subset $G\subset\BBR^N$, then $u(x,t)\geq \gf_\infty(t)$.
\es
\Proof By assumption, there exists a sequence $\{t_n\}$ decreasing to $0$ such that 
	\bel{J-2}
	\lim_{n\to \infty}\myint{G}{}u(x,t_n)dx=\infty.
  	\ee
If $(\ref{J-1})$ holds, we can construct a decreasing sequence of open subsets $G_k\subset G$ such that $\overline {G_k}\subset G_{k-1}$, diam$(G_k)=\ge_k\to 0$ when $k\to\infty$, and
	\bel{J-3}
	\lim_{n\to \infty}\myint{G_k}{}u(x,t_n)dx=\infty \forevery k\in\BBN.
	\ee
Furthermore there exists a unique $a\in\cap_k G_k$. We set 
$$\myint{G_k}{}u(x,t_n)dx=M_{n,k}.
$$
Since ${\dsps \lim_{n\to\infty}M_{n,k}=\infty}$, we claim that for any $m>0$ and any $k$, there exists $n=n(k)\in\BBN$ such that 
    \bel{J-4}\myint{G_k}{}u(x,t_{n(k)})dx\geq m.  \ee
By induction, we define $n(1)$ as the smallest integer $n$ such that $M_{n,1}\geq m$. This is always possible. Then we define  $n(2)$ as the smallest integer larger than $n(1)$  such that $M_{n,2}\geq m$. By induction, $n(k)$ is the smallest integer $n$ larger than $n(k-1)$ such that $M_{n,k}\geq m$. Next, for any $k$, there exists $\ell=\ell (k)$ such that
    \bel{J-5}\myint{G_k}{}\inf\{u(x,t_{n(k)});\ell\}dx= m   \ee
and we set 
	$$\hat U_{k}(x)=\inf\{u(x,t_{n(k)});\ell\}\chi_{_{G_k}}(x). $$ 
Let $\hat u_k=u$ be the unique bounded solution of 
	\bel{J-6}\left\{\BA {lll}
	\prt_tu-\Gd_p u+f(u)&=0\qq&\text{in }Q_{\infty}\\\phantom{,,,v+f(v)}
	u(.,0)&=\hat U_{k}\qq&\text{in }\BBR^N.
	\EA\right. \ee
Since $\hat u_k(x,0)\leq u(x,t_{n(k)})$, we derive  
	\bel{J-7}  u(x,t+t_{n(k)}) \geq \hat u_k(x,t)\qq\forall (x,t)\in Q_\infty.   \ee
When $k\to\infty$, $\hat U_{k}\rightarrow m\gd_a$, thus $\hat u_k\to u_{m\gd_a}$ by \rth{stab}. Therefore $u\geq u_{m\gd_a}$. Since $m$ is arbitrary and $ u_{m\gd_a}\to\phi_\infty$ when $m\to\infty$, it follows that 
$u\geq \phi_\infty$.\qeda
  
\blemma {Jinfinite}Assume $p>2$, $(\ref{CFS})$ and $J=\infty$ are satisfied, and ${\dsps \lim_{k \to \ity}u_k}=\ity$. There exists no positive solution $u$ of $(\ref{A1})$ in $Q_{\infty}$ which satisfies $(\ref {J-1})$ for some bounded open subset $G\subset\BBR^N$.
\es
\Proof If we assume that such a $u$ exists, we proceed as in the proof of the previous lemma. Since \rth{stab} holds, we derive that $u\geq u_{m\gd_a}$ for any $m$. Since ${\dsps \lim_{m\to\infty}u_{m\gd_a}(x,t)=\infty}$ for all $(x,t)\in Q_\infty$, we are led to a contradiction.\qeda  \medskip

Thanks to these results, we can characterize the initial trace of positive solutions of $(\ref{A1})$ when the Keller-Osserman condition does not hold. \medskip

\noindent{\bf Proof of \rth{IS2}. }{\it (i)} If $\CS(u)\ne \emptyset$, there exists $y \in \CS(u)$ and an open neighborhood $G$ of $y$ such that $(\ref{J-1})$ holds. By \rlemma{Jfinite}, $u\geq \phi_\infty$ and the initial trace of $u$ is the Borel measure $\gn_\infty$. Otherwise, $\CR(u)=\BBR^N$ and $Tr_{_{\BBR^N}}(u) \in \GTM_+(\BBR^N)$. \smallskip

\noindent{\it (ii)} Using the argument as in \rth{tr+J-fin} and because of \rlemma{Jinfinite}, $\CS(u)=\emptyset$. Therefore $\CR(u)=\BBR^N$ and $Tr_{_{\BBR^N}}(u) \in \GTM_+(\BBR^N)$. \qeda

\bcor{IS3} Assume $p>2$. If $f$ is convex and satisfies $(\ref{CFS})$, $J<\infty$ and $K=\infty$, there exist infinitely many different positive solutions $u$ of $(\ref{A1})$ such that $tr_{\BBR^N}(u)=\gn_\infty$. 
\es
\Proof Let $b>0$ be fixed. Since $f$ is increasing, $(x,t)\mapsto U(x,t)=w_b(x)+\phi_{\infty(t)}$ is a supersolution for  $(\ref{A1})$. Let $V(x,t)=\max\{w_b(x),\phi_{\infty}(t)\}$ then $V$, $f(V)$ and $|\nabla V|^p$ are locally integrable in $Q_T$; actually $V$ is locally Lipschitz continuous. Let
$\ge>0$ and $\gr_\ge$ be a smooth approximation defined by
$$
\gr_\ge(r)=\left\{\BA {lll}
0\qquad&\text{if }\;& r<0\\
\frac{r^2}{2\ge}&\text{if }\;&0<r<\ge\\
r-\frac{\ge}{2}\qquad&\text{if }\;&r>\ge
\EA\right.
$$
We set $V_\ge(x,t)=\phi_\infty(t)+\gr_\ge[w_b(x)-\phi_\infty(t)]$. Then
$$\BA {l}
\prt_tV_\ge-\Gd_pV_\ge+f(V_\ge)=\phi'_\infty\left(1-\gr_\ge'[w_b-\phi_\infty]\right)
-\left(\gr_\ge'[w_b-\phi_\infty]\right)^{p-1}\Gd_pw_b
\\\phantom{\prt_tV_\ge-\Gd_pV_\ge+f(V_\ge)}
-(p-1)\left(\gr_\ge'[w_b-\phi_\infty]\right)^{p-2}\gr_\ge''[w_b-\phi_\infty]|\nabla w_b|^p+f(V_\ge)
\\\phantom{\prt_tV_\ge-\Gd_pV_\ge+f(V_\ge)}
\leq f(V_\ge)-\left(1-\gr_\ge'[w_b-\phi_\infty]\right)f(\phi_\infty)-\left(\gr_\ge'[w_b-\phi_\infty]\right)^{p-1}f(w_b)
\EA$$

If $\gf\in C^\infty_c(Q_T)$ is nonnegative, then
$$\myint{}{}\myint{Q_T}{}\left(-V_\ge\prt_t\gf+|\nabla V_\ge |^{p-2}\nabla V_\ge.\nabla \gf+f(V_\ge)\right)dx\,dt\leq o(1)
$$
Letting $\ge \to 0$ implies
$$\myint{}{}\myint{Q_T}{}\left(-V\prt_t\gf+|\nabla V |^{p-2}\nabla V.\nabla \gf+f(V)\right)dx\,dt\leq 0.
$$
 Thus $V$ is a subsolution, smaller than $U$. Therefore there exists a solution $u_b$ such that
 $V\leq u\leq U$. This implies that $tr_{\BBR^N}(u_b)=\gn_\infty$.  If $b'>b$ we construct $u_{b'}$ with 
 $tr_{\BBR^N}(u_{b'})=\gn_\infty$ and $\lim_{t\to\infty}(u_{b'}(0,t)-u_b(0,t))>0$.
\qeda

\end{document}